\documentclass[12pt, reqno]{amsart} 
\usepackage[normalem]{ulem}
\usepackage[a4paper, margin=2.5cm, top=3.5cm, bottom=3.5cm]{geometry}
\usepackage[colorlinks=true,
linkcolor=webgreen,
filecolor=webbrown,
citecolor=webgreen]{hyperref}
\usepackage{array}
\newcolumntype{C}[1]{>{\centering\arraybackslash}m{#1}}
\newenvironment{longlist}
 {
  \par\nobreak\vspace{\abovedisplayskip}%
  \begingroup\lccode`~=`, \lowercase{\endgroup\def~}{,\hfil\penalty0\hfilneg}%
  \catcode`,=\active
  \let\originaltimes\times
  \renewcommand{\times}{\originaltimes{}\discretionary{}{\kern1em}{}}%
  \setlength{\parindent}{0pt}%
  \hangindent=3.4em \hangafter=1
 }
 {\par\vspace{\belowdisplayskip}}
\usepackage{amssymb}
\usepackage{amsmath}
\usepackage{amsthm}
\usepackage{mathtools}
\DeclarePairedDelimiter\ceil{\lceil}{\rceil}
\DeclarePairedDelimiter\floor{\lfloor}{\rfloor}
\usepackage[usenames]{color}
\usepackage{graphicx}
\usepackage[dvipsnames]{xcolor}

\definecolor{webgreen}{rgb}{0,.5,0}
\definecolor{webbrown}{rgb}{.6,0,0}
\definecolor{red}{rgb}{1,0,0}

\usepackage[normalem]{ulem}

\usepackage{color}
\theoremstyle{definition}
\newtheorem{definition}{Definition}[section]
\newcommand{\seqnum}[1]{\href{http://oeis.org/#1}{\underline{#1}}}
\begin{document}
\theoremstyle{plain}
\newtheorem{theorem}[definition]{Theorem}
\newtheorem{corollary}[definition]{Corollary}
\newtheorem{lemma}[definition]{Lemma}
\newtheorem{proposition}[definition]{Proposition}
\newtheorem{fact}[definition]{Fact}
\theoremstyle{definition}
\newtheorem{example}[definition]{Example}
\newtheorem{conjecture}[definition]{Conjecture}
\theoremstyle{bolditalic}
\newtheorem{remark}[definition]{Remark}
\newtheorem{open question}[definition]{Open question}
\title{On $k$-layered numbers}
\author[F. Jokar]{Farid Jokar}
\email{jokar.farid@gmail.com}

\date{}

\maketitle

\begin{abstract}
A positive integer $n$ is said to be $k$-layered if its divisors can be partitioned into $k$ sets with equal sum. In this paper, we start the systematic study of these class of numbers. In particular, we state some algorithms to find some even $k$-layered numbers $n$ such that $2^{\alpha}n$ is a $k$-layered number for every positive integer $\alpha$. We also find the smallest $k$-layered number for $1\leq k\leq 8$. Furthermore, we study when $n!$ is a $3$-layered and when is a $4$-layered number. Moreover, we classify all $4$-layered numbers of the form $n=p^{\alpha}q^{\beta}rt$, where $\alpha$, $1\leq \beta\leq 3$, $p$, $q$, $r$, and $t$ are two positive integers and four primes, respectively. In addition, in this paper, some other results concerning these numbers and their relationship with $k$-multiperfect numbers, near-perfect numbers, and superabundant numbers are discussed. Also, we find an upper bound for the differences  of two consecutive $k$-layered numbers for every positive integer $1\leq k\leq 5$. Finally, by assuming the smallest $k$-layered number, we find an upper bound for the difference of two consecutive $k$-layered numbers. 
\end{abstract}
\section{Introduction}
A perfect number is a positive integer that is equal to the sum of its proper positive divisors. In 2003, the idea of Zumkeller numbers, which is  a generalization of perfect numbers, was first introduced by Zumkeller in the Online Encyclopedia of Integer Sequences (OEIS) \cite{OEIS} (\seqnum{A083207}).
\begin{definition}
A positive integer $n$ is said to be Zumkeller if the set of all positive divisors of $n$ can be partitioned into two subsets  with the same sum. A Zumkeller partition for a Zumkeller number $n$ is
a partition $\{A_1, A_2\}$ of the set of all positive divisors of $n$ such that $A_1$ and $A_2$ sums to the same value. 
\end{definition}
Shortly thereafter, Clark et al. \cite{clark} announced several results and conjectures related to Zumkeller numbers. Peng and Bhaskara Rao \cite{zumkeller} found some other results about Zumkeller numbers; they also studied the relationship between practical numbers (to be defined later in the paper) and Zumkeller numbers. They also settled a conjecture from \cite{clark}. Moreover, they made substantial contributions towards the second conjecture of Clark et al. \cite{clark}. Recently, the author studied a generalization of Zumkeller numbers called $k$-layered numbers \cite{difference} which is the principal object of study in this paper.

\begin{definition}
A positive integer $n$ is said to be  $k$-layered if the set of all positive divisors of $n$ can be partitioned into $k$  subsets with the same sum. A $k$-layered partition for a $k$-layered number $n$ is
a partition $\{A_1, A_2, \dots, A_k\}$ of the set of all positive divisors of $n$ such that for every $1\leq i\leq k $, each of $A_i$ sums to the same value. 
\end{definition}

In \cite{Mahanta}, the authors found some new results related to Zumkeller numbers and $k$-layered numbers. It is clear from the definition of $k$-layered numbers that all Zumkeller numbers are $2$-layered numbers. Another generalization of perfect numbers which we will consider are multiperfect numbers or $k$-multiperfect number and near-perfect numbers. A positive integer $n$ is called a multiperfect number or a $k$-multiperfect number if $n$ satisfies the equation $\sigma(n)=kn$, for some positive integer $k$. It is easy to check that every perfect number is Zumkeller; this immediately gives us a further motivation to define $k$-layered numbers as all multiperfect numbers will be $k$-layered numbers.

This paper is arranged as follows: in Sections \ref{basic} and \ref{almost section}, we recall and generalize some results related to $k$-layered numbers, taking as basis the results already proved for Zumkeller numbers. In Section \ref{the smallest}, we find the smallest $k$-layered number for $1\leq k\leq 8$. In Section \ref{some more}, we state some algorithms to find some even $k$-layered numbers $n$ such that $2^{\alpha}n$ is $k$-layered for every positive integer $\alpha$.  In addition, we prove that $n!$ is a $3$-layered number if and only if $n\geq 5$ and $n\neq10$. Moreover, we  prove that  $n!$ is a $4$-layered number if and only if $n\geq 9$. In Section \ref{multiperfect number section}, we prove that all known $3$-multiperfect numbers are $3$-layered. Furthermore, by assuming  Conjecture \ref{semiperfect}, we prove that all known $4$-multiperfect numbers and all known $5$-multiperfect numbers are $4$-layered and $5$-layered, respectively. In addition, we state some algorithms to find some more $k$-multiperfect numbers that are $k$-layered. In Section \ref{near perfect number section}, we find some relationship between Near-perfect numbers and Zumkeller numbers.  In Section \ref{difference section}, we find an upper bound for  the differences of two consecutive $k$-layered number for $1\leq k\leq 5$.
 Lastly, by assuming the smallest $k$-layered number, we find an upper bound for the differences of two consecutive $k$-layered numbers. There is a generous sprinkling of open questions and conjectures in this paper, which the author hopes subsequent mathematicians will take up.

\section{Some preliminary results on $k$-layered numbers}\label{basic}
Let $n$ be a positive integer and $\sigma(n)$ denote the sum of positive divisors of $n$. We recall that the abundancy of $n$ is defined to be $I(n)=\dfrac{\sigma(n)}{n}$. Also, A number $n$ is said to be  abundant, perfect, and deficient if $I(n)>2$, $I(n)=2$, and $I(n)<2$, respectively. Table \ref{tab:1} lists the smallest $k$-layered number for every positive integer $2 \leq k \leq 8$ (for more details, see Section \ref{the smallest} of this paper).

\begin{table}[h!]
	\centering
    \label{tab:table2}
    \begin{tabular}{|C{2cm}|c|C{2.5cm}|}
    \hline
    \multicolumn{1}{|C{2cm}|}{ $k$ }
    & Smallest $k$-layered number & Number of positive divisors\\ \hline
    $1$ & $1$ & $1$ \\ \hline
    $2$ & $6$ & $4$ \\ \hline
    $3$ & $120$ & $16$ \\ \hline
    $4$ & $27720$ & $96$ \\ \hline
    $5$ &$147026880$ & $896$\\ \hline
    $6$ & $130429015516800$ & $18432$ \\ \hline
    $7$ & $1970992304700453905270400$ &$1474560$ \\ \hline
    $8$ & $1897544233056092162003806758651798777216000$& $1245708288
    $ \\ \hline
   \end{tabular}
\caption{Smallest $k$-layered number for $1\leq k\leq 8$.}
\label{tab:1}
\end{table}
\begin{remark}
 According to the number of divisors of the smallest 8-layered number $n$, it seems that it is not easy to find an 8-layered partition for $n$. The author applied an advanced algorithm to find an 8-layered partition for $n$.
\end{remark}
In the remainder of this paper, if $A$ is a set of all positive integers, then we define $S(A)$ to be the sum of the integers in $A$. We do not prove the results in this section because the proofs are exactly similar to the proofs of the corresponding results for Zumkeller numbers in \cite{zumkeller}. The following proposition gives a necessary and sufficient condition for a positive integer $n$ to be $k$-layered.
 \begin{proposition}\label{3layered}
A number $n$ is $k$-layered if and only if we can find $k-1$ disjoint subsets $A_1, A_2,\dots, A_{k-1}$ of positive divisors of $n$ so that for every $1\leq i\leq k-1$, $A_i$ sums to the $\dfrac{\sigma(n)}{k}$.
\end{proposition}

Proposition $2$ from \cite{zumkeller} gives some necessary conditions for a Zumkeller number. We can generalize this proposition for $k$-layered number.

\begin{proposition}\label{clear}
If $n$ is a $k$-layered number, then the following is true:

$\mathrm{(i)}\ $ $ k\mid\sigma (n)$;
    
$\mathrm{(ii)}\ $ If $k$ is even, then the prime factorization of $n$ must include at least one odd prime to an odd power;

$\mathrm{(iii)}\ $ $ \sigma (n)\geq kn$; this concludes that $I(n)\geq k$.
\end{proposition}
\begin{proof}
The proof is identical to the proof of Proposition $2$ of  \cite{zumkeller}.
\end{proof}

\begin{remark}\label{grows}
Let $m$ be the smallest $9$-layered number. We can see that the number 
$$n=4368924363354820808981210203132513655327781713900627249499856876120704000$$
 is the smallest number such that $I(n)\geq 9$ (\seqnum{A023199}). Therefore, by Proposition \ref{clear}, $m\geq n $.  
\end{remark}
Although it is an open question that for every positive integer $k$, there exists a $k$-layered numbers, but Table \ref{tab:1} and Remark \ref{grows} lead us to the following conjecture.
\begin{conjecture}\label{expon}
Size of the smallest $k$-layered number grows exponentially with respect to $k$.
\end{conjecture}
We close this section with two propositions which are a generalization of Corollary 5 and Proposition 6 of \cite{zumkeller}, respectively.

\begin{proposition} \label{nw}
If the integer $n$ is $k$-layered and $w$ is relatively prime to $n$, then $nw$ is a $k$-layered number.
\end{proposition}
\begin{proof}
 The proof is identical to the proof of  Corollary $5$ of \cite{zumkeller}.
\end{proof}

\begin{proposition}
Let $n$ be a $k$-layered number and $p_1^{a_1}p_2^{a_2}\cdots p_m^{a_m}$ be the prime factorization of $n$. Then for any positive integers $l_1, l_2, \ldots, l_m$, the number
\[
p_1^{a_1+l_1(a_1+1)}p_2^{a_2+l_2(a_2+1)}\cdots p_1^{a_m+l_m(a_m+1)}
\] is $k$-layered.
\end{proposition}
\begin{proof}
The proof of the above result is similar to the proof of Proposition 6 in \cite{zumkeller}.
\end{proof}

\section{Almost practical numbers and $k$-layered numbers}\label{almost section}

In this section we discuss the connection between $k$-layered numbers and almost practical numbers, which are a generalization of practical numbers. The results discussed here will be a generalization of the results connecting Zumkeller numbers and practical numbers that were found by Peng and Bhaskara Rao \cite{zumkeller}.
 
First, we recall the definition of practical numbers.
\begin{definition}
A positive integer $n$ is said to be a practical number if every positive integer less than $n$ can be represented as a sum of distinct positive divisors of $n$.
\end{definition}
Now we recall two results about practical numbers.
    
\begin{proposition}[Corollary 1, \cite{practical}]\label{prop5}
A positive integer $n$ with the prime factorization $p_1^{t_1} p_2^{t_2} \cdots p_m^{t_m}$ and $p_1 < p_2 < \dots < p_m$ is a practical number if and only if $p_1 = 2$ and $p_{i+1} \leq \sigma(p_1^{t_1} \cdots p_i^{t_i}) + 1$ for $1 \leq i \leq m - 1$.
\end{proposition}
\begin{proposition}[Proposition 8, \cite{zumkeller}]\label{prop6}
A positive integer $n$ is a practical number if and only if every integer less than or equal to $\sigma(n)$ can be written as a sum of distinct divisors of $n$.
\end{proposition}
We now define almost practical numbers, which are a generalization of practical numbers. These numbers were first studied by Stewart \cite{practical}.
\begin{definition}
A positive integer $n$ is called an almost practical number if all of the numbers $j$ such that
$2<j<\sigma(n) -2$ can be written as a sum of distinct divisors of $n$. 
\end{definition}

\begin{remark}\label{every}
It is clear that every practical number is an almost practical number.
\end{remark}
Now we recall some results related to almost practical numbers from \cite{practical}.
\begin{proposition}[Theorem 3, \cite{practical}] \label{either}
Let $n \neq 3$ be an odd positive integer and $1=d_1<d_2<\dots<d_k=n$ are the divisors of $n$. We also define $\sigma_i=d_1+d_2+\cdots +d_i$. Then, $n$ is an almost practical number if and only if $d_2=3, d_3=5$ and for $i\geq 3$, at least one of the following is true:
\begin{enumerate}
    \item $ d_{i+1}\leq \sigma_i-2$ and $ d_{i+1}\neq \sigma_i-4 $,
    \item $ d_{i+1}= \sigma_i-4$ and $ d_{i+2}=\sigma_i-2$.
\end{enumerate}
\end{proposition}

\begin{remark}\label{remark almost} 
If $n=p_1^{\alpha_1} p_2^{\alpha_2} \dots p_m^{\alpha_m}$ is an odd almost practical number where $p_1 < p_2 < \dots < p_m$ are all prime factors of $n$, then by Proposition \ref{either}, it is clear that $p_3=7$.
\end{remark}

\begin{theorem}[Theorem 4, \cite{practical}]\label{almost} 
Let $n \neq 3$ be an odd almost practical number and $p$ be an odd prime, then $pn$ is an almost practical number if and only if $2p \leq \sigma(n)-2$ and $2p \neq \sigma(n)-4$.
\end{theorem}

\begin{remark}
By Propositions \ref{either} and \ref{prop5}, if $n$ is an odd almost practical number, then $2n$ is a practical number.
\end{remark}

\begin{proposition}[Corollary 4, \cite{practical}]\label{almost divisors}
Let $n \neq 3$ be an odd almost practical number and $p$ be a prime dividing $n$, then $pn$ is an almost practical number.
\end{proposition}
\begin{example} \label{f2}
By  Proposition \ref{either}, one can check directly that $3^3\times 5 \times 7$ is an almost practical number. Therefore, by Proposition \ref{almost divisors}, if  $n=3^{\alpha_1}\times 5^{\alpha_2}\times 7^{\alpha_3} $ such that $\alpha_1\geq 3$, $\alpha_2$, and $\alpha_3$ are positive integers, then the number $n$ is an almost practical number.
\end{example}
\begin{example}\label{f22}
 Let $k\geq 4$ be a positive integer. Also, let $p_i$ denotes $i$th prime number. By Bertrand postulate, for every integer $i \geq 4$, we have
$$2p_i< 4p_{i-1}\leq  \sigma(p_2p_3\cdots p_{i-1})-4.$$
Thus, by Example \ref{f2} and Theorem \ref{almost}, if $m=p_2^{\alpha_1} p_3^{\alpha_2}\cdots p_k^{\alpha_{k-1}}$ such that $\alpha_1\geq3$, $\alpha_2, \ldots, \alpha_{k-1}$ are positive integers, then $m$ is an almost practical number.
\end{example}
Now we are going to investigate the relation between almost practical numbers and Zumkeller numbers.
The following proposition is a generalization of Proposition 10 of \cite{zumkeller}.
\begin{proposition}\label{17}
Let $n \neq 3$ be an almost practical number. Then, $n$ is Zumkeller if and only if $\sigma(n)$ is even.
\end{proposition}
\begin{proof} 
The proof is similar to the proof of Proposition 10 of \cite{zumkeller}.
\end{proof}
\begin{example}\label{ex5}
 Let $k\geq 4$ be a positive integer, and also let $p_i$ denotes $i$th prime number. Now let $m=p_2^{\alpha_1} p_3^{\alpha_2}\cdots p_k^{\alpha_{k-1}}$ such that $\alpha_1\geq 3$, $\alpha_2, \ldots$ and $\alpha_{k-1}$ are positive integers. Thus, by Example \ref{f22} and Proposition \ref{17}, $m$ is a Zumkeller number if and only if we can find an integer $1\leq i \leq k-1$ such that $\alpha_i$ is an odd number. In other words, the number $m$ is Zumkeller if and only if $m$ is not square. 
\end{example}
Now we state a proposition which is a generalization of Proposition $13$ of \cite{zumkeller}.
\begin{proposition}\label{nonzumkeller}
 Let $\alpha$ be a positive integer. Also, let $n$ be non-$k$-layered number  and $p$ be a prime number with $\gcd(n, p)=1$. If $np^{\alpha}$ is a $k$-layered number, then $p\leq \sigma(n)$. Moreover, if $q$ is prime factor of $k$ such that  $q\nmid\sigma(n)$, then $q\mid\sigma(p^{\alpha})$. 
\end{proposition}
\begin{proof}
The proof is identical to the proof of  Proposition $13$ of \cite{zumkeller}.
\end{proof}
In  Proposition 20 of \cite{zumkeller}, the authors by a simple method proved that every odd Zumkeller number possesses at least three distinct prime factors. On the other hand, we know that for every positive integer $t$, the number $2^{t}$ is a deficient number. Therefore, by Propositions \ref{clear} and \ref{nonzumkeller}, if $p$ is a prime number such that $\sigma(2^{t})<p$, then the number $2^{t}p^{\alpha}$ fails to be a Zumkeller number for every positive integer $\alpha$. Also, by Proposition \ref{17}, if $p\leq \sigma(2^{t})$ and $\alpha$ is an odd number, then the number $2^{t}p^{\alpha}$ is Zumkeller. This was also proved in \cite{Mahanta} (see Theorem 2.6) by a different method. Thus, we have the following corollary.
\begin{corollary}\label{2p}
Let $p_1$ and $p_2$ be prime numbers. Also, let $\alpha_1$ and $\alpha_2$ be positive integers. The number $p_1^{\alpha_1}p_2^{\alpha_2}$ is Zumkeller if and only if  $p_1=2$, $p_2\leq \sigma(2^{\alpha_1})$ and $\alpha_2$ is odd.
\end{corollary}
\begin{remark}\label{Pankaj}
Let $p$ be an odd prime number. Also. let $\alpha$ and $\beta$ be positive integers. One can check that the number $2^{\alpha}p^{\beta}$ is abundant if and only if $p\leq \sigma(2^{\alpha})$ (for more details, see the proof of Theorem 2.6 of \cite{Mahanta}). Thus, if $n$ is an abundant number that possesses two factors, then $n$ is Zumkeller if and only if $n$ satisfies the conditions $\mathrm{(i)}\ $  and $\mathrm{(iii)}\ $ of Proposition \ref{clear}.
\end{remark}
\section{The smallest $k$-layered number for $1 \leq k\leq 8$}\label{the smallest}
    We start this section with two definitions.
     \begin{definition}
   A positive integer $n$ is said to be $k$-abundant if $I(n)\geq k$.
   \end{definition}
   \begin{definition}
   A positive integer $n$ is said to be good $k$-abundant if $I(n)\geq k$ and $k\mid \sigma(n)$.
   \end{definition}
   \begin{remark}
    It is clear that a positive integer $n$ is 2-abundant if and only if $n$ is either perfect or abundant.
   \end{remark}
    Let $i$ be a positive integer. In remainder of this section, we define $a_i$, $g_i$, and $k_i$ as the smallest $i$-abundant number, the smallest good $i$-abundant number, and the smallest $i$-layered number, respectively. In \cite{OEIS}, Pegg listed $a_i$ for $1\leq i\leq 13$ (\seqnum{A023199}). Let $i$ be a positive integer. The following simple proposition states a relationship between $a_i$, $g_i$, and $k_i$.
    \begin{proposition}\label{bound}
     For positive integers $i$, $a_i\leq g_i\leq k_i$.
    \end{proposition}
   \begin{proof}
   By Proposition \ref{clear}, it is obvious.
   \end{proof}
    We recall that a positive integer $n$ is said to be superabundant if $I(n)>I(m)$ for all positive integers $m<n$. Now we recall a well-know theorem about superabundant numbers of \cite{AE}.
 \begin{theorem}\label{tsup}
 If $n$ is a superabundant number, then there exists a positive integer $k$ such that $n=p_1^{\alpha_1}p_2^{\alpha_2}\dots p_{k}^{\alpha_{k}}$, where $p_i$ is $i$th prime number and $\alpha_1\geq \alpha_2\geq \dots \geq \alpha_{k}\geq 1$. Also, $a_k=1$ except when $n$ is $4$ or $36$.
\end{theorem}
Thus, by previous theorem, every superabundant is almost practical (see Example \ref{f22}). Therefore, as a direct consequence of the previous theorem and Proposition \ref{17}, we have the following corollary.
\begin{corollary}
A positive integer $n$ is superabundant if and only if one of the following holds:

$\mathrm{(i)}\ $ $n$ is Zumkeller;

$\mathrm{(ii)}\ $ $n=4$  or $n=36$.
\end{corollary}
A list of all superabundant numbers less than $10^{1200}$ can be found in \cite{superabundant}. The following proposition states one of the special properties of superabundant numbers.
   \begin{proposition}\label{pitsu}
 Let $m_1<m_2$ be two consecutive superabundant numbers. For every positive  integer $t< m_2$, $I(t)\leq I(m_1)$, and the equality holds when $t=m_1$.
 \end{proposition}
 \begin{proof}
 If $t\leq m_1,$ then $I(t)\leq I(m_1)$. So, we assume $m_1<t_1<m_2,$ such that $I(t_1)>I(m_1).$ Since $m_1$ and $m_2$ are consecutive superabundant numbers, therefore $t_1$ is not a superabundant number. So, there exists $t_2<t_1$ such that $I(t_2)>I(t_1)$. Hence $I(t_2)>I(t_1)>I(m_1),$ which implies $m_1<t_2<t_1$. We once again conclude that there exists an integer $t_3$ where $m_1<t_3<t_2<t_1$ such that $I(m_1)<I(t_1)<I(t_2)<I(t_3)$. Thus, for every positive integer $r\geq3$, by this algorithm, inductively, we can find distinct positive integers $t_1, t_2, \dots,$ and $t_r$ such that $m_1<t_r<t_{r-1}<\dots<t_1<m_2$ and $I(m_2)> I(t_r)> I(t_{r-1})> \dots>I(t_1)>I(m_1)$; this contradicts the finiteness of the set $\{a: a\in \mathbb{N}, m_1<a<m_2\}$.
 \end{proof}
 In remainder of this section, we define $s_i$ as the $i$th superabundant number. As a direct consequence of Propositions \ref{bound} and \ref{pitsu}, we have the following corollary.
  \begin{corollary}\label{algorithm}
For positive integer $i$, there exists a unique positive integer $t_i$ such that satisfies the following two conditions:
 
 $\mathrm{(i)}\ $ $I(s_{\ell})\geq i$ for integers $\ell\geq t_i$,

 $\mathrm{(ii)}\ $ $I(s_{j})<i$ for integers $j<t_i$.
 
  Also, if a positive integer $t$ satisfies conditions $\mathrm{(i)}\ $ and  $\mathrm{(ii)}\ $ of the proposition, then $s_t=a_i$. Furthermore, $a_i=g_i$ or $s_{t+1}\leq g_i$.
  \end{corollary}
Let $n$ be a positive integer. In the remainder of this paper, we define $D_n$ as the set all positive divisors of $n$. Proposition \ref{bound} and  corollary \ref{algorithm} motivate an algorithm to find $a_i$, $g_i$, and $k_i$ for positive integers $i$.
 \begin{example}\label{exgood}
 By the list of the smallest $k$-abundant numbers in \cite{OEIS} (\seqnum{A023199}) and the list of superabundant numbers \cite{superabundant},  $a_5=s_{43}=122522400=2^{5}\times3^{2}\times5^{2}\times7\times11\times13\times17$ and $s_{44}=147026880=2^6 \times 3^3 \times 5 \times 11 \times 13 \times 17$. One can check directly that $5\nmid \sigma(s_{43})$. Therefore, $a_5\neq g_5$. Also, $5\mid \sigma(s_{44})$. Thus, by  corollary \ref{algorithm}, $g_5=s_{44}$. Now we prove that $g_5$ is $5$-layered. We define
\begin{longlist}
$A_{1} = \{$ \begin{math} 3^{} \end{math}, \ \ 
\begin{math} 2^{6} \end{math}, \ \ 
\begin{math} 2^{6} \times 3^{} \end{math}, \ \ 
\begin{math} 3^{2} \times 5^{} \times 17^{} \end{math}, \ \ 
\begin{math} 2^{6} \times 11^{} \times 13^{} \end{math}, \ \ 
\begin{math} 2^{6} \times 3^{} \times 5^{} \times 7^{} \times 17^{} \end{math}, \ \ 
\begin{math} 2^{6} \times 3^{3} \times 11^{} \times 17^{} \end{math}, \ \ 
\begin{math} 2^{6} \times 3^{3} \times 5^{} \times 7^{} \times 11^{} \times 13^{} \times 17^{} \end{math}$\},$
\end{longlist}

\begin{longlist}
$A_{2} = \{$ \begin{math} 2^{2} \times 17^{} \end{math}, \ \ 
\begin{math} 2^{6} \times 7^{} \times 13^{} \end{math}, \ \ 
\begin{math} 2^{2} \times 3^{3} \times 5^{} \times 7^{} \times 11^{} \times 17^{} \end{math}, \ \ 
\begin{math} 2^{6} \times 3^{} \times 5^{} \times 7^{} \times 13^{} \times 17^{} \end{math}, \ \ 
\begin{math} 2^{5} \times 3^{3} \times 5^{} \times 11^{} \times 13^{} \times 17^{} \end{math}, \ \ 
\begin{math} 2^{5} \times 3^{2} \times 5^{} \times 7^{} \times 11^{} \times 13^{} \times 17^{} \end{math}, \ \ 
\begin{math} 2^{4} \times 3^{3} \times 5^{} \times 7^{} \times 11^{} \times 13^{} \times 17^{} \end{math}, \ \ 
\begin{math} 2^{5} \times 3^{3} \times 5^{} \times 7^{} \times 11^{} \times 13^{} \times 17^{} \end{math}$\},$
\end{longlist}
\begin{longlist}
$A_{3} = \{$ \begin{math} 5^{} \end{math}, \ \ 
\begin{math} 2^{4} \end{math}, \ \ 
\begin{math} 3^{2} \times 5^{} \times 7^{} \end{math}, \ \ 
\begin{math} 2^{6} \times 5^{} \end{math}, \ \ 
\begin{math} 2^{6} \times 3^{2} \end{math}, \ \ 
\begin{math} 2^{3} \times 3^{} \times 5^{} \times 7^{} \end{math}, \ \ 
\begin{math} 2^{6} \times 3^{3} \end{math}, \ \ 
\begin{math} 2^{4} \times 3^{} \times 7^{} \times 11^{} \times 17^{} \end{math}, \ \ 
\begin{math} 2^{6} \times 3^{3} \times 5^{} \times 13^{} \end{math}, \ \ 
\begin{math} 2^{6} \times 7^{} \times 11^{} \times 13^{} \times 17^{} \end{math}, \ \ 
\begin{math} 2^{4} \times 3^{2} \times 7^{} \times 11^{} \times 13^{} \times 17^{} \end{math}, \ \ 
\begin{math} 2^{6} \times 3^{2} \times 5^{} \times 11^{} \times 13^{} \times 17^{} \end{math}, \ \ 
\begin{math} 2^{6} \times 3^{3} \times 5^{} \times 7^{} \times 11^{} \times 17^{} \end{math}, \ \ 
\begin{math} 2^{4} \times 3^{2} \times 5^{} \times 7^{} \times 11^{} \times 13^{} \times 17^{} \end{math}, \ \ 
\begin{math} 2^{6} \times 3^{3} \times 5^{} \times 7^{} \times 13^{} \times 17^{} \end{math}, \ \ 
\begin{math} 2^{5} \times 3^{3} \times 7^{} \times 11^{} \times 13^{} \times 17^{} \end{math}, \ \ 
\begin{math} 2^{6} \times 3^{} \times 5^{} \times 7^{} \times 11^{} \times 13^{} \times 17^{} \end{math}, \ \ 
\begin{math} 2^{3} \times 3^{3} \times 5^{} \times 7^{} \times 11^{} \times 13^{} \times 17^{} \end{math}, \ \ 
\begin{math} 2^{6} \times 3^{3} \times 5^{} \times 11^{} \times 13^{} \times 17^{} \end{math}, \ \ 
\begin{math} 2^{6} \times 3^{3} \times 7^{} \times 11^{} \times 13^{} \times 17^{} \end{math}$\},$
\end{longlist}
\begin{longlist}
$A_{4} = \{$ \begin{math} 2^{2} \end{math}, \ \ 
\begin{math} 2^{4} \times 3^{} \times 11^{} \end{math}, \ \ 
\begin{math} 2^{} \times 3^{3} \times 5^{} \times 7^{} \times 17^{} \end{math}, \ \ 
\begin{math} 2^{6} \times 3^{3} \times 5^{} \times 13^{} \times 17^{} \end{math}, \ \ 
\begin{math} 2^{3} \times 3^{3} \times 7^{} \times 11^{} \times 13^{} \times 17^{} \end{math}, \ \ 
\begin{math} 2^{6} \times 3^{2} \times 5^{} \times 7^{} \times 11^{} \times 17^{} \end{math}, \ \ 
\begin{math} 2^{4} \times 3^{} \times 5^{} \times 7^{} \times 11^{} \times 13^{} \times 17^{} \end{math}, \ \ 
\begin{math} 2^{6} \times 3^{3} \times 11^{} \times 13^{} \times 17^{} \end{math}, \ \ 
\begin{math} 2^{5} \times 3^{3} \times 5^{} \times 7^{} \times 11^{} \times 13^{} \end{math}, \ \ 
\begin{math} 2^{6} \times 3^{2} \times 5^{} \times 7^{} \times 13^{} \times 17^{} \end{math}, \ \ 
\begin{math} 2^{} \times 3^{3} \times 5^{} \times 7^{} \times 11^{} \times 13^{} \times 17^{} \end{math}, \ \ 
\begin{math} 2^{5} \times 3^{2} \times 7^{} \times 11^{} \times 13^{} \times 17^{} \end{math}, \ \ 
\begin{math} 2^{4} \times 3^{3} \times 5^{} \times 11^{} \times 13^{} \times 17^{} \end{math}, \ \ 
\begin{math} 2^{6} \times 5^{} \times 7^{} \times 11^{} \times 13^{} \times 17^{} \end{math}, \ \ 
\begin{math} 2^{5} \times 3^{3} \times 5^{} \times 7^{} \times 11^{} \times 17^{} \end{math}, \ \ 
\begin{math} 2^{3} \times 3^{2} \times 5^{} \times 7^{} \times 11^{} \times 13^{} \times 17^{} \end{math}, \ \ 
\begin{math} 2^{5} \times 3^{3} \times 5^{} \times 7^{} \times 13^{} \times 17^{} \end{math}, \ \ 
\begin{math} 2^{4} \times 3^{3} \times 7^{} \times 11^{} \times 13^{} \times 17^{} \end{math}, \ \ 
\begin{math} 2^{5} \times 3^{} \times 5^{} \times 7^{} \times 11^{} \times 13^{} \times 17^{} \end{math}, \ \ 
\begin{math} 2^{6} \times 3^{3} \times 5^{} \times 7^{} \times 11^{} \times 13^{} \end{math}, \ \ 
\begin{math} 2^{2} \times 3^{3} \times 5^{} \times 7^{} \times 11^{} \times 13^{} \times 17^{} \end{math}, \ \ 
\begin{math} 2^{6} \times 3^{2} \times 5^{} \times 7^{} \times 11^{} \times 13^{} \times 17^{} \end{math}$\},$
\end{longlist}
\noindent and $A_5=D_{g_5}\setminus (A_1\cup A_2\cup A_3\cup A_4)$, then one can check directly that $\{A_1, A_2, A_3, A_4, A_5\}$ is a $5$-layered partition for $g_5$. Thus, by Proposition \ref{bound}, $k_5=g_5$.
 \end{example} 
 \begin{remark}\label{goodcon}
One can check directly that $3\mid\mid a_2$, $5\mid\mid a_3$, $7\mid\mid a_4$, $29\mid\mid a_6$, $13\mid\mid a_7$, $23\mid\mid a_8$, $53\mid\mid a_9$, $269\mid\mid a_{10}$, $461\mid\mid a_{11}$, $47\mid\mid a_{12}$, and $1481\mid\mid a_{13}$. This concludes that  if $i$ is a positive integer such that $i\neq 1 ,5$ and $2\leq i \leq 13$, then there exists a prime number $q_i$ such that $q_i\mid\mid a_i$ and $q_i\equiv i-1$ (mod $i$). Thus, if $i$ is a positive integer such that $i\neq 5$ and $1\leq i \leq 13$, then $i\mid \sigma(a_i)$. In other words, for integers $1\leq i \leq 13$ except $i=5$, $g_i=a_i$.
 \end{remark}
 Remark \ref{goodcon} and Example \ref{exgood} lead us to the following open question.
  \begin{open question}
For positive integers $k\neq5$, dose $a_i$ equal $g_i$?
\end{open question}
 In the following, we prove that $g_i=k_i$ for $1\leq i \leq 8$.
 \begin{proposition}
  $g_i=k_i$ for  integers $1\leq i \leq 8$.
 \end{proposition}
 \begin{proof}
  \textbf{Case 1}: \ $k=1$. It is clear.
  
 \textbf{Case 2}: \ $k=2$. We know that $a_2=g_2=6$ is the smallest Zumkeller number. 

\textbf{Case 3}: \ $k=3$. By the list of the smallest $k$-abundant numbers in \cite{OEIS} (\seqnum{A023199}), the list of superabundant numbers in \cite{superabundant}, and Corollary \ref{algorithm}, $s_{10}=120=a_3=g_3$. Now if we define $A_1=\{$\begin{math}20\end{math}, \ \ \begin{math}40\end{math}, \ \ \begin{math}60\end{math}$\}$, $A_2=\{$\begin{math}1\end{math}, \ \ \begin{math}2\end{math}, \ \ \begin{math}3\end{math}, \ \ \begin{math}4\end{math}, \ \ \begin{math}5\end{math}, \ \ \begin{math}6\end{math}, \ \ \begin{math}8\end{math}, \ \ \begin{math}10\end{math}, \ \ \begin{math}12\end{math}, \ \ \begin{math}15\end{math}, \ \ \begin{math}24\end{math}, \ \ \begin{math}30\end{math}$\}$, and $A_3=\{ $\begin{math} 120 \end{math}$\}$, then $\{A_1, A_2, A_3\}$ is a $3$-layered partition for $g_3$. Thus, by Proposition \ref{bound}, $k_3=g_3$.

\textbf{Case 4}: \ $k=4$. By the same method that we used in Case 3, We can see $s_{23}=a_4=g_4=27720=2^{3}\times3^2\times5\times7\times11$. Now we define
\begin{longlist}
$A_{1} = \{$ \begin{math} 2^{3} \times 3^{2} \times 5^{} \end{math}, \ \ 
\begin{math} 2^{3} \times 3^{2} \times 5^{} \times 7^{} \times 11^{} \end{math}$\},$
\end{longlist}
\begin{longlist}
$A_{2} = \{$ \begin{math} 2^{} \end{math}, \ \ 
\begin{math} 2^{3} \times 3^{} \end{math}, \ \ 
\begin{math} 2^{2} \times 3^{} \times 5^{} \end{math}, \ \ 
\begin{math} 2^{} \times 5^{} \times 7^{} \end{math}, \ \ 
\begin{math} 2^{2} \times 3^{} \times 5^{} \times 11^{} \end{math}, \ \ 
\begin{math} 2^{} \times 5^{} \times 7^{} \times 11^{} \end{math}, \ \ 
\begin{math} 2^{2} \times 3^{2} \times 5^{} \times 7^{} \end{math}, \ \ 
\begin{math} 2^{} \times 3^{2} \times 7^{} \times 11^{} \end{math}, \ \ 
\begin{math} 2^{2} \times 5^{} \times 7^{} \times 11^{} \end{math}, \ \ 
\begin{math} 2^{3} \times 3^{} \times 7^{} \times 11^{} \end{math}, \ \ 
\begin{math} 2^{2} \times 3^{2} \times 5^{} \times 11^{} \end{math}, \ \ 
\begin{math} 2^{2} \times 3^{} \times 5^{} \times 71^{} \end{math}, \ \ 
\begin{math} 2^{2} \times 3^{2} \times 5^{} \times 7^{} \times 11^{} \end{math}$\},$
\end{longlist}
\begin{longlist}
$A_{3} = \{$ \begin{math} 1 \end{math}, \ \ 
\begin{math} 2^{} \times 3^{2} \end{math}, \ \ 
\begin{math} 2^{} \times 3^{} \times 5^{} \times 7^{} \times 11^{} \end{math}, \ \ 
\begin{math} 2^{2} \times 3^{2} \times 7^{} \times 11^{} \end{math}, \ \ 
\begin{math} 2^{3} \times 5^{} \times 7^{} \times 11^{} \end{math}, \ \ 
\begin{math} 3^{2} \times 5^{} \times 7^{} \times 11^{} \end{math}, \ \ 
\begin{math} 2^{3} \times 3^{2} \times 5^{} \times 11^{} \end{math}, \ \ 
\begin{math} 2^{3} \times 3^{2} \times 7^{} \times 11^{} \end{math}, \ \ 
\begin{math} 2^{} \times 3^{2} \times 5^{} \times 7^{} \times 11^{} \end{math}$\},$
\end{longlist}
\noindent and $A_{4} = D_{g_4}\setminus{(A_1 \cup A_2 \cup A_3)}$, then one can check directly that $\{A_1, A_2, A_3, A_4\}$ is a $4$-layered partition for $g_4$.

\textbf{Case 5}:  $k=5$. By Example \ref{exgood}, it is clear.

\textbf{Case 6}:  $6\leq k\leq 8$. In \cite{small}, the author stated $6$-layered partition, $7$-layered partition, and $8$-layered partition for $g_6=130429015516800=2^{7} \times 3^{3} \times 5^{2} \times 7^{2} \times 11^{} \times 13^{} \times 17^{} \times 19^{} \times 23^{} \times 29^{}, g_7= 1970992304700453905270400=2^{7} \times 3^{4} \times 5^{2} \times 7^{2} \times 11^{} \times 13^{} \times 17^{} \times 19^{} \times 23^{} \times 29^{} \times 31^{} \times 37^{} \times 41^{} \times 43^{} \times 47^{} \times 53^{}$, and $g_8=1897544233056092162003806758651798777216000=2^{10} \times 3^{5} \times 5^{3} \times 7^{2} \times 11^{2} \times 13^{} \times 17^{} \times 19^{} \times 23^{} \times 29^{} \times 31^{} \times 37^{} \times 41^{} \times 43^{} \times 47^{} \times 53^{} \times 59^{} \times 61^{} \times 67^{} \times 71^{} \times 73^{} \times 79^{} \times 83^{} \times 89^{}$, respectively.
\end{proof}
The previous proposition leads us to the following conjecture.
\begin{conjecture}
 $g_i=k_i$ for every positive integer $i$.
\end{conjecture}
 \section{Some more results about $k$-layered numbers}\label{some more}
In \cite{zumkeller}, the authors if a number $n$ is Zumkeller, then the number $2^{\alpha}n$ is half-Zumkeller \footnote{A positive integer $n$ is said to be a half-Zumkeller number if the proper positive divisors of n can be partitioned into two disjoint subsets of equal sum. A half-Zumkeller partition for a half-Zumkeller number $n$ is a partition $\{A_1 , A_2\}$ of the set of proper positive divisors of $n$ such that $S(A_1)=S(A_2)$.} for every positive integer $\alpha$. By the same method, in \cite{Mahanta}, the authors proved that if a number $n$ is Zumkeller, then the number $2^{\alpha}n$ is Zumkeller for every positive integer $\alpha$. This does not generalize to all $k$-layered numbers. For instance, the number $m=2^{5}\times3\times7$ is a $3$-layered number because if we define $A_1=\{2^{5}\times3\times7\}$, $A_2=\{2^{4}\times3\times7, 2^3\times3\times7, 2^2\times3\times7, 2\times3\times7, 3\times7, 2\times7, 2\times3, 1\}$, then one can check directly that $A_1=A_2=\dfrac{\sigma(m)}{3}$. This concludes that the number $m$ is $3$-layered. However, the number $2m$ fails to be $3$-layered since $\sigma(2m)$ is not divisible by $3$. By Proposition \ref{nw}, we know that if $n$ is an odd $k$-layered number, then for positive integer $\alpha$ the number $2^{\alpha}n$ is a $k$-layered number. This leads us to the following definition.
\begin{definition}
An even $k$-layered number $n$ is said to be a good $k$-layered number if the number $2^{\alpha}n$ is $k$-layered for every positive integer $\alpha$.
\end{definition}
The following proposition states a necessary condition for the $k$-layered number $n$ to be a good $k$-layered.
\begin{proposition}
Let $\alpha$ be a positive integer, and let $m$ be an odd positive integer. If $n=2^{\alpha}m$ is a good $k$-layered number, then $k\mid\sigma(m)$.
\end{proposition}
\begin{proof}

Let $n$ be a good $k$-layered number such that $k\nmid\sigma(m).$ Then there exists at least one factor $p^{\gamma}$ of $k$ such that $p^{\gamma}\nmid\sigma(m)$ and $p^{\gamma}\mid\sigma(2^{\alpha}),$ where $p$ is a prime number and $\gamma$ is a positive integer. If $p=2$, then this is a contradiction. If $p$ is odd, then this implies that $p^{\gamma}\nmid2^{\alpha+1}.$ This concludes that $p^{\gamma}\nmid\sigma(2^{\alpha+1})$. Thus, $p^{\gamma}\nmid\sigma(2n),$ which is a contradiction.
\end{proof}

In the following, we state two crucial theorems that help us to find a subset of good $k$-layered numbers.
For better understanding, first, we state a proposition, which is a special case of one of these two theorems.
\begin{proposition}\label{is 3layered}
Let $n$ be an odd positive integer and $\alpha$ be a positive integer. Let $m=2^{\alpha}n$ be a $3$-layered number with $3$-layered partition $\{A_1, A_2, A_3\}$. Now let $A'_1$ be a subset of the set of all positive divisors of $n$ such that $A'_1$ sums to $\dfrac{2\sigma(n)}{3}$. If $A'=\{2^{\alpha}d: d \in A'_1\}$ and $A'\subset A_1 \cup A_2$, then $m$ is a good $3$-layered number.
\end{proposition}
\begin{proof}
We define $M_1=A'\cap A_1$ and $M_2=A' \cap A_2$. Now we define
 \[M'_1=\{2d: d \in M_1\}, M'_2=\{2d: d\in M_2\}, M_3=\{2^{\alpha+1}d: d\in (D_n\setminus A'_1) \}.\]
 
We can see that $$S((A_1\setminus M_1)\cup M'_1)=S(A_1 \cup M_1).$$  Thus, we see
 \begin{multline*}
      S(A_1)+\dfrac{2^{\alpha+1}\sigma(n)}{3}=S((A_1\setminus M_1)\cup M'_1\cup M_2)=S(A_1\cup M_1\cup M_2)=S(A_2\cup M_2\cup M_1)\\=S((A_2\setminus M_2)\cup M'_2\cup M_1) =S(A_2)+\dfrac{2^{\alpha+1}\sigma(n)}{3}.
 \end{multline*}

Therefore, by definition of $A'$, one can check directly that $$\{(A_1\setminus M_1)\cup M'_1\cup M_2,(A_2\setminus M_2)\cup M'_2\cup M_1, A_3 \cup M_3\}$$ is a $3$-layered partition for $2m$. Thus, by the same method, inductively, we can prove that $\ell=2^{t}n$ is a $3$-layered number for every integer $t\geq \alpha$.
\end{proof}
By Proposition \ref{is 3layered},  we have the following corollary.
\begin{corollary}\label{3layered2}
If $n\neq3$ is an odd almost practical number such that $6\mid\sigma(n)$, then the number $2n$ is a good $3$-layered number.
\end{corollary}
\begin{proof}
First, we prove that $2n$ is a $3$-layered  number. By Proposition \ref{17},  $n$ is a Zumkeller number; this concludes that $D_n$  can be partitioned into two subsets $A_1$ and $A_2$ such that each of them sums to $\dfrac{\sigma(n)}{2}$. Now we define $B_1=\{2d: d\in A_1\}$ and $B_2=\{2d: d\in A_2\}$. We know that $n$ is an odd number. Therefore, for every integer $a\in B_1 \cup B_2$, $a\notin D_n$. Thus, $\{D_n, B_1, B_2\}$ is a $3$-layered partition for $2n$. Also, we know that $n$ is an almost practical number such that $6\mid \sigma(n)$ and $2\neq\dfrac{2\sigma(n)}{3}\neq \sigma(n)-2$; this concludes that there exists $A'\subset D_n$ so that $A'$ sums to $\dfrac{2\sigma(n)}{3}$. Now let $A''=\{2d: d\in A'\}$. We know that $A'' \subset B_1 \cup B_2$. Then, by Proposition  \ref{is 3layered}, $2n$ is a good $3$-layered number.
\end{proof}

Now we present a proposition generalizing Proposition $\ref{is 3layered}$ for $4$-layered numbers.
 \begin{proposition}\label{is 4layered}
 Let $n$ be an odd Zumkeller number with Zumkeller partition $\{A'_1, A'_2\}$. For positive integer $\alpha$, we define $A''_1=\{2^{\alpha}d: d\in A'_1\}$ and $A''_2=\{2^{\alpha}d: d\in A'_2\}$. Now let  $m=2^{\alpha}n$ be a $4$-layered number with $4$-layered partition $\{A_1, A_2, A_3, A_4\}$. If $A''_1\subset (A_1\cup A_2)$ and $A''_2\subset (A_3\cup A_4)$, then the number $m$ is good $4$-layered.
 \end{proposition}
 \begin{proof}
 We define\\
 
 $M_1=A''_1\cap A_1, M_2=A''_1 \cap A_2$\\
 
 $M_3=A''_2\cap A_3, M_4=A''_2 \cap A_4$\\
 
Now for every integer $1\leq i \leq 4$, we define $M'_i=\{2d: d \in M_i\}.$ One can check directly that $$\{(A_1\setminus{M_1}) \cup M'_1\cup M_2, (A_2\setminus M_2)\cup M'_2\cup M_1,(A_3\setminus M_3)\cup M'_3\cup M_4, (A_4\setminus{M_4})\cup M'_4\cup M_3)\}$$ is a $4$-layered partition for $m$. Thus, by applying the same method like before, inductively, we can prove that $\ell=2^{t}n$ is a $4$-layered number for every integer $t\geq \alpha$.
 \end{proof}
 
The following two results are a generalization of Propositions \ref{is 3layered} and \ref{is 4layered}, respectively. 
\begin{theorem}\label{good1}
Let $k\geq 3$ and $n$ be two odd positive integers such that $k\mid \sigma(n)$. Let $A'_1,A'_2,\dots, A'_{\frac{k-1}{2}}$ be disjoint subsets the set of all positive divisors of $n$ so that $A'_i$ sums to $\dfrac{2\sigma(n)}{k}$ for every integer $1\leq i\leq \dfrac{k-1}{2}$. Let $\alpha$ be a positive integer and $A''_i=\{2^{\alpha}d: d\in A'_i\}$ for every integer $1\leq i\leq \dfrac{k-1}{2}$. Now let $m=2^{\alpha}n$ be a $k$-layered number with $k$-layered partition $\{A_1, A_2, \dots, A_k\}$ such that $A''_i\subset A_{2i-1}\cup A_{2i}$ for every integer $1\leq i\leq \dfrac{k-1}{2}$. Then, $m$ is a good $k$-layered number.
\end{theorem}
 \begin{proof}
For integers $1\leq i \leq k-1$, we define

$$M_{i}=A''_{\floor*{\frac{i+1}{2}}} \cap A_i, M'_i= \{2d: d \in M_{i}\},$$ and $$M=\{2^{\alpha+1}d : d \in (D_n\setminus A'_1 \cup A'_2 \cup \dots , A'_{\frac{K-1}{2}})\}.$$

Now we define 

\begin{displaymath}
B_i=\begin{cases}
   (A_i\setminus{M_{i})\cup M'_{i}\cup M_{i+1} }, & \text {if $i$ is odd;} \\
   (A_i\setminus{M_{i})\cup M'_{i}\cup M_{i-1}}, & \text {if $i$ is even.}
\end{cases}
\end{displaymath}
One can check directly that $\{B_1, B_2,\dots, B_{k-1}, A_k\cup M\}$ is a $k$-layered partition for $2^{\alpha +1}n$. Also, by applying the previous method, inductively, we can prove that $\ell = 2^{t}n$ is $k$-layered for every integer $t\geq \alpha$.
\end{proof}

\begin{theorem}\label{generalization of good even}
Let $k\geq4$ be an even positive integer. Let  the odd positive integer $n$ be a $\dfrac{k}{2}$-layered number with $\dfrac{k}{2}$-layered partition $\{A'_1,A'_2,\dots, A'_{\frac{k}{2}}\}$. Let $\alpha$ be a positive integer  and $A''_i=\{2^{\alpha}d: d\in A'_i\}$ for every integer $1\leq i\leq \dfrac{k}{2}$. Also, let $m=2^{\alpha}n$ be a $k$-layered number with $k$-layered partition $\{A_1, A_2, \dots, A_k\}$ so that $A''_i\subset A_{2i-1}\cup A_{2i}$ for every integer $1\leq i\leq \dfrac{k}{2}$. Then, the number $m$ is a good $k$-layered number.
\end{theorem}
\begin{proof}
For integers $1 \leq i \leq k$, we define

$$M_i=A''_{\floor*{ \frac{i+1}{2}}}\cap A_i, M'_i=\{2d: d\in M_i\}.$$

Now for every $1\leq i \leq k $, we define

\begin{displaymath}
B_i=
\begin{cases}
   (A_i\setminus{M_{i})\cup M'_{i}\cup M_{i+1} }, & \text {if $i$ is odd;} \\
   (A_i\setminus{M_{i})\cup M'_{i}\cup M_{i-1}}, & \text {if $i$ is even.}
\end{cases}
\end{displaymath}
It is easy to check that the set $\{B_1, B_2, \dots , B_k\}$ is $k$-layered partition for $2^{\alpha+1}n$. Also, by applying the same method, inductively, we can prove the number $2^{t}n$ is $k$-layered for every $t\geq \alpha$.
\end{proof}
As a consequence of Theorems \ref{good1} and \ref{generalization of good even}, we have the following corollary.
\begin{corollary}\label{good corollary}
For integers $2\leq k\leq 5$, the smallest $k$-layered number is a good $k$-layered number. 
\end{corollary}
\begin{proof}
We investigate the five cases for $k$.

\textbf{Case 1}: \ $k=2$. This was already proved in Proposition $4.12$ of \cite{Mahanta}.

\textbf{Case 2}: \ $k=3$. Let $\{A_1, A_2, A_3\}$ be the $3$-layered partition for the smallest $3$-layered number, $120$, which was defined in Section \ref{the smallest}. Now let $n=120=2^3m$. We know that $\dfrac{2\sigma(m)}{3}=16$. If we define $A'_1=\{1,15\}$, then it is obvious that $A'=\{2^3, 2^3\times 15\} \subset A_2 \cup A_3 $. Thus, by Proposition \ref{is 3layered}, $n$ is a good $3$-layered number. 

\textbf{Case 3}: \ $k=4$. Let $\{A_1, A_2, A_3, A_4\}$ be the $4$-layered partition for the smallest $4$-layered $n=27720$, which was defined in Section \ref{the smallest}. Now let $D'$ be the set of all positive divisors of $m=3465=3^2\times5\times7\times11$. We define

\begin{longlist}
$A'_{1} = \{$ \begin{math} 3^{} \end{math}, \ \ 
\begin{math} 3^{2} \times 5^{} \end{math}, \ \ 
\begin{math} 3^{} \times 7^{} \times 11^{} \end{math}, \ \ 
\begin{math} 3^{2} \times 5^{} \times 7^{} \times 11^{} \end{math}$\},$
\end{longlist}
\noindent and $A'_2=D'\setminus{A'_1}$.

One can check directly that for every integer $1\leq i \leq 2$, $S(A'_i)=\dfrac{\sigma(m)}{2}$. Therefore, $\{A'_1, A'_2\}$ is a Zumkeller partition for $m$. We know that $n=2^3m$. Also, if we define $A''_1=\{2^3d: d\in A'_1\}$ and $A''_2=\{2^3d: d\in A'_2\}$, then one can check directly that $A''_1\subset A_1 \cup A_2$ and $A''_2 \subset A_3 \cup A_4$. Thus, by Proposition \ref{is 4layered}, $n$ is a good $4$-layered number.

\textbf{Case 4}: \ $k=5$. Let $\{A_1, A_2, A_3, A_4, A_5\}$ be the $5$-layered partition for the smallest $5$-layered number $n=147026880$, which was defined in Section \ref{the smallest}. Let $n=2^6m$. We define
\begin{longlist}
$A'_{1} = \{$ \begin{math} 1 \end{math}, \ \ 
\begin{math} 3^{} \end{math}, \ \ 
\begin{math} 11^{} \times 13^{} \end{math}, \ \ 
\begin{math} 3^{} \times 5^{} \times 7^{} \times 17^{} \end{math}, \ \ 
\begin{math} 3^{} \times 5^{} \times 7^{} \times 13^{} \times 17^{} \end{math}, \ \ 
\begin{math} 2^{6} \times 3^{3} \times 5^{} \times 7^{} \times 11^{} \times 13^{} \times 17^{} \end{math}$\},$
\end{longlist}
\noindent and
\begin{longlist}
$A'_{2} = \{$ \begin{math} 5^{} \end{math}, \ \ 
\begin{math} 3^{2} \end{math}, \ \ 
\begin{math} 3^{3} \end{math}, \ \ 
\begin{math} 3^{3} \times 5^{} \times 13^{} \end{math}, \ \ 
\begin{math} 7^{} \times 11^{} \times 13^{} \times 17^{} \end{math}, \ \ 
\begin{math} 3^{2} \times 5^{} \times 11^{} \times 13^{} \times 17^{} \end{math}, \ \ 
\begin{math} 3^{3} \times 5^{} \times 7^{} \times 11^{} \times 17^{} \end{math}, \ \ 
\begin{math} 3^{} \times 5^{} \times 7^{} \times 11^{} \times 13^{} \times 17^{} \end{math}, \ \ 
\begin{math} 3^{3} \times 5^{} \times 11^{} \times 13^{} \times 17^{} \end{math}, \ \ 
\begin{math} 3^{3} \times 7^{} \times 11^{} \times 13^{} \times 17^{} \end{math}$\}.$
\end{longlist}

One can check directly that $S(A'_1)=S(A'_2)=\dfrac{2\sigma(m)}{5}$. We define $A''_i=\{2^6d : d\in A'_i\}$ for $i=1, 2$. We know $A'_1\subset (A_1 \cup A_2)$ and $A'_2\subset (A_3 \cup A'_4)$. Thus, by Theorem \ref{good1}, $n$ is a good $5$-layered.
\end{proof}

Corollary \ref{good corollary}, leads us naturally to the following open question.
\begin{open question}
Is the smallest $k$-layered number a good $k$-layered number for every positive integer $k\neq 1$?
\end{open question}
It appears that the methods used in this paper are not feasible to attack this problem at the moment.

Now we state a proposition that can be used to construct some $k$-layered numbers.
\begin{proposition}\label{kzum}
Let $k$ and $k'$ be positive integers. Now let $m$ be a $k$-layered number and  $n$ be a $k'$-layered number. If $\gcd(m, n) =1 $, then $mn$ is a $kk'$-layered number.
\end{proposition}
\begin{proof}
 Let $\{A_1 ,A_2,\dots,A_k\}$  be a $k$-layered partition for $m$, and $\{B_1, B_2, \dots, B_{k'}\}$ be a $k'$-layered partition for n. One can check directly that $\{A_iB_j: 1\leq i \leq k , 1\leq j \leq k'\}$ is a $kk'$-layered partition for $mn$.
\end{proof}
It is clear that Proposition \ref{kzum} can be generalized.
\begin {corollary}
Let $\ell_1, \ell_2, \dots, \ell_r$ be positive integers such that for every integer $1\leq i \leq r$, $m_i$ is a $\ell_i$-layered number. If for every distinct integers $1\leq i\neq j \leq r$, $\gcd(m_i, m_j)=1$, then $m_1m_2\dots m_r$ is a $\ell_1\ell_2\dots \ell_r$-layered number.
 \end{corollary}
The following example shows the power of Proposition \ref{kzum} for finding a subset of the set of $4$-layered numbers.
\begin{example}
 By Example \ref{ex5}, $n_1=3^{\alpha_1}\times 5^{\alpha_2}\times 7^{\alpha_3}$, in which $\alpha_1\geq3$, $\alpha_2$, $\alpha_3$ are positive integers, and at least one of the exponents of its factors is odd, is a Zumkeller number. Let $t$ be a positive integer and $p$ be a prime number such that $p\leq 2^{t+1}-1$ and $\gcd(p, n_1)=1$. By Corollary \ref{2p}, for every odd number $\alpha_4$, the number $n_2=2^t\times p^{\alpha_4}$ is Zumkeller. Therefore, by  Proposition \ref{kzum}, $n=n_1n_2$ is a $4$-layered number.
 \end{example}
In \cite{zumkeller}, the authors proved that the number $n!$ is Zumkeller for integers $n\geq3$ (we can also prove it in a slightly different way, see Remark \ref{faczum}). In the following, we prove that the number $n!$ is $3$-layered for every integer $n\geq 5$ and $n\neq10$, and also the number $n!$ is $4$-layered for every integer $n\geq 9$. Before that, we  recall a theorem which was proved by Breusch \cite{B}; this theorem is a generalization of Bertrand's postulate .
\begin{theorem}[Breusch \cite{B}]\label{Breusch}
 For every integer $n \geq 7$, there are primes of the form $3k+1$ and $3k+2$ between $n$ and $2n$.
 \end{theorem}
\begin{theorem}\label{f3}
Let $n\geq11$ be an integer. Then, the  number $n!$ possesses prime factorization $p^{\alpha_1}_1p^{\alpha_2}_2\dots p^{\alpha_k}_k$ such that $2=p_1<p_2<\dots < p_k$ and $\alpha_{k-1}=\alpha_{k}=1$. In addition, $p_k \leq 2^{\alpha_1}$ and there exists a prime number $q_n$ such that $q_n\mid\mid n!$ and $q_n\equiv 2$ (mod 3).
 \end{theorem}
 \begin{proof}
 If $11\leq n\leq 13$, then one can check directly that $n!$ satisfies the theorem. Now Let $n\geq 14$. Also, let $p'$ be a prime number. We know that $2p'\leq n$ if and only if $p'^{2}\mid n!$. Thus, if $p$ is the largest prime number such that $p^2\mid n!$, then by definition of $n$, it is clear that $p\geq 7$. Therefore, by Theorem \ref{Breusch}, there exist at least two distinct prime numbers $p'_1$ and $p'_2$ such that $p<p'_1 , p'_2 <2p$ and for at least one $i\in\{1,2\},$ $p'_i\equiv 2$ (mod $3$). Also, by definition of $p$, $p'_1\mid\mid n!$ and $p'_2\mid\mid n!$. Furthermore, if $\nu_{2}(n!)$   denotes the exponent of the largest power of $2$ that divides $n!$, then by Legendre's formula,  we have
 $$p_k \leq n < 2^{\floor*{\frac{n}{2}}} <2^{\nu_{2}(n!)}.$$
 \end{proof}
Now we state and prove one of the key theorems of this section.
  \begin{theorem}\label{factorial 3}
  The number $n!$ is a $3$-layered number if and only if $n\geq 5$ and $n\neq10$.
  \end{theorem}
  \begin{proof}
    On can check directly that for every integer $1\leq n \leq 4$, $I(n!)<3$, and $3\nmid\sigma(10!)$. So, by Proposition \ref{clear}, $n!$ fails to be $3$-layered for every integer $1\leq n \leq 4$ and $n=10$.
    
 Now we  investigate the five cases for $n$:
    
\textbf{Case 1}: $n=5$. By Corollary \ref{good corollary}, $n!=120$ is a $3$-layered number.

\textbf{Case 2}: $n=6$. Let $A_1=\{2\times 3,2^4\times 5,2^4\times 3^2\times 5\}$ and $A_2=\{2^3,2^3\times 3,2^4\times 3^2,2^3\times 3\times 5,2^4\times 3\times 5,2\times 3^2\times 5,2^2\times 3^2\times 5\}$. One can check directly that $S(A_1)=S(A_2)=\dfrac{\sigma(n!)}{3}.$ Thus, by Proposition \ref{3layered}, the number $n!=2^4\times 3^2\times 5$ is a $3$-layered number.

\textbf{Case 3}: $n=7$. We proved that the number $6!=2^{4}\times3^2\times5$ is a $3$-layered number. Thus, by Proposition \ref{nw}, the number $7!=2^{4}\times3^2\times5\times7$ is a $3$-layered number.

\textbf{Case 4}: $n=8$. We define

$A_{1} = \{$ \begin{math} 2^{3} \times 3^{} \times 5^{} \end{math}, \ \ 
\begin{math} 2^{3} \times 3^{2} \times 5^{} \times 7^{} \end{math}, \ \ 
\begin{math} 2^{5} \times 3^{2} \times 5^{} \times 7^{} \end{math}, \ \ 
\begin{math} 2^{7} \times 3^{2} \times 5^{} \times 7^{} \end{math}$\}$  

$A_{2} = \{$ \begin{math} 2^{3} \end{math}, \ \ 
\begin{math} 2^{3} \times 3^{} \times 7^{} \end{math}, \ \ 
\begin{math} 2^{7} \times 5^{} \times 7^{} \end{math}, \ \ 
\begin{math} 2^{6} \times 3^{} \times 5^{} \times 7^{} \end{math}, \ \ 
\begin{math} 2^{7} \times 3^{2} \times 7^{} \end{math}, \ \ 
\begin{math} 2^{7} \times 3^{} \times 5^{} \times 7^{} \end{math}, \ \ 
\begin{math} 2^{6} \times 3^{2} \times 5^{} \times 7^{} \end{math}$\}$\\
One can check directly that $S(A_1)=S(A_2)=\dfrac{\sigma(n!)}{3}$. Thus, by Proposition \ref{3layered}, the number $n!=2^{7}\times3^2\times5\times7$ is a $3$-layered number.

\textbf{Case 5}: $n\geq9$. Let $n=9$. Then, $n!=2^{7}m$, where $m=3^{4}\times5\times7$. By Example \ref{f2}, $m$ is an odd almost practical number such that $6\mid \sigma(m)$.  Therefore, by Corollary \ref{3layered2}, $n!$ is a $3$-layered number. Now let $n>9$ and $n!=2^{\alpha_1}m'$ with $\gcd(2^{\alpha_1}, m')=1$ such that $m'$ is an odd almost practical number. Also, let $p^{\alpha_2}_2p^{\alpha_3}_3\cdots p^{\alpha_k}_k$ be the prime factorization of $m'$. We want to prove that $(n+1)!=2^{\beta_1}m''$ with $\gcd(2^{\beta_1}, m'')=1$ such that $m''$ is an odd almost practical number.
  First, let $n+1$ be a composite number. Then, if $p$ is a prime number such that $p\mid n+1$, then there exists a positive integer $d\neq 1$ such that $n+1=pd$. This concludes that $p<n$. Therefore, $p\mid n!$. Thus, $(n+1)!=2^{\beta_1}p^{\beta_2}_2\cdots p^{\beta_k}_k$, where $\beta_i\geq \alpha_i$ for every integer $1\leq i \leq k$. Therefore, by Proposition \ref{almost divisors}, $m''=p^{\beta_2}_2\cdots p^{\beta_k}_k$ is an odd almost practical number. Moreover, if $n+1=q$ is a prime number, then by  Bertrand postulate, we have
  $$2q< 4p_k< \sigma(p_2p_3\cdots p_{k})-4.$$
  By Theorem \ref{almost}, we once again conclude that $m'' $ is an odd almost practical number.
  Thus, for every integer $n\geq 9$, $n!=2^{\alpha_1}t$, where $t$ is an odd almost practical number. Also, by Theorem \ref{f3}, there exists a prime number $q_n$ such that $q_n\mid\mid n!$ and $q_n\equiv 2$ (mod $3$) for all $n\geq11$. Therefore, there exists a positive integer $t'$ such that $t=q_nt'$ and $\gcd(q_n, t')=1$; this concludes that $6\mid\sigma(q_n)\sigma(t')=\sigma(t)$. 
  Thus, by Corollary \ref{3layered2}, for every integer $n\geq9$ except $n=10$, n! is a $3$-layered number.
  \end{proof}

  \begin{remark}\label{38f}
  In fact, in case 5, by Corollary \ref{3layered2}, we proved that the number $n!$ is a good $3$-layered number for $n\geq 11$. One can check that $n!$ is a good $3$-layered number for $5\leq n\leq 9$, we omit the details here. For instance, by Corollary \ref{good corollary}, the number $5!$ is a good $3$-layered number. Thus, we can say that the number $n!$ is $3$-layered if and only if $n!$ is a good $3$-layered number.
  \end{remark}
 
 \begin{remark}\label{39f}
  Let $n\geq9$ and $n\neq10$. Using the same arguments of Case 5 above, we get the number $\dfrac{n!}{2^t}$ is a $3$-layered number for every integer $0\leq t \leq \nu_{2}(n!)-1$. Thus, by the previous remark, if $r=\nu_{2}(n!)-1$, then the number $\dfrac{n!}{2^r}$ is a good $3$-layered number for every integer $n\geq9$ and $n\neq10$.
  \end{remark}

  \begin{theorem}\label{factorial 4}
 The number $n!$ is $4$-layered number if and only if $n\geq 9$.
  \end{theorem}
  \begin{proof}
  On can check directly that for every integer $1\leq n \leq 8$, $I(n!)<4$. Then, by Proposition \ref{clear}, $n!$ fails to be $4$-layered for every integer $1\leq n \leq 8$. Now we want to prove that for every integer $n\geq 9$ the number $n!$ is $4$-layered.
We investigate the three cases for $n$:\\
\textbf{Case 1}: \ $n=9$. We define

$A_{1} = \{$ \begin{math} 2^{2} \times 3^{} \end{math}, \ \ 
\begin{math} 2^{3} \times 3^{4} \end{math}, \ \ 
\begin{math} 2^{6} \times 3^{} \times 5^{} \times 7^{} \end{math}, \ \ 
\begin{math} 2^{7} \times 3^{4} \times 5^{} \times 7^{} \end{math}$\}$

$A_{2} = \{$ \begin{math} 2^{2} \end{math}, \ \ 
\begin{math} 2^{7} \times 7^{} \end{math}, \ \ 
\begin{math} 2^{4} \times 3^{4} \times 5^{} \end{math}, \ \ 
\begin{math} 2^{6} \times 3^{3} \times 5^{} \times 7^{} \end{math}, \ \ 
\begin{math} 2^{7} \times 3^{3} \times 5^{} \times 7^{} \end{math}, \ \ 
\begin{math} 2^{6} \times 3^{4} \times 5^{} \times 7^{} \end{math}$\}$

$A_{3} = \{$ \begin{math} 2^{2} \times 3^{2} \end{math}, \ \ 
\begin{math} 2^{6} \times 3^{2} \times 5^{} \end{math}, \ \ 
\begin{math} 2^{5} \times 3^{3} \times 5^{} \times 7^{} \end{math}, \ \ 
\begin{math} 2^{6} \times 3^{4} \times 7^{} \end{math}, \ \ 
\begin{math} 2^{7} \times 3^{2} \times 5^{} \times 7^{} \end{math}, \ \ 
\begin{math} 2^{4} \times 3^{4} \times 5^{} \times 7^{} \end{math}, \ \ 
\begin{math} 2^{7} \times 3^{4} \times 5^{} \end{math}, \ \ 
\begin{math} 2^{7} \times 3^{4} \times 7^{} \end{math}, \ \ 
\begin{math} 2^{5} \times 3^{4} \times 5^{} \times 7^{} \end{math}$\}$

 One can check directly that for every integer $1\leq i \leq 3 $, $S(A_i)=\dfrac{\sigma(n!)}{4}$. Thus, by Proposition \ref{3layered}, the number $n!=362880=2^{7}\times3^{4}\times5\times7$ is a $4$-layered number.
 
 \textbf{Case 2}: \ $n=10$. We define
 
 $A_{1} = \{$\begin{math} 2^{} \end{math}, \begin{math} 2^{5} \times 3^{} \end{math}, \ \ \begin{math} 2^{4} \times 3^{3} \times 7^{} \end{math}, \ \
\begin{math} 2^{7} \times 3^{2} \times 5^{2} \times 7^{} \end{math}, \ \ \begin{math} 2^{8} \times 3^{4} \times 5^{2} \times 7^{} \end{math}$\}$

$A_{2}=\{$\begin{math} 1 \end{math}, \ \ \begin{math}
 2^{4} \end{math}, \ \ 
\begin{math} 3^{2} \times 5^{2} \end{math}, \ \ 
\begin{math} 2^{6} \times 3^{2} \times 5^{} \end{math}, \ \ 
\begin{math} 2^{8} \times 3^{2} \times 5^{} \times 7^{} \end{math}, \ \ 
\begin{math} 2^{8} \times 3^{4} \times 5^{} \times 7^{} \end{math}, \ \
\begin{math} 2^{8} \times 3^{3} \times 5^{2} \times 7^{} \end{math}, \ \ 
\begin{math} 2^{7} \times 3^{4} \times 5^{2} \times 7^{} \end{math}$\}$

$A_{3}=\{$\begin{math} 2^{} \times 3^{2} \end{math}, \ \
\begin{math} 2^{5} \times 7^{} \end{math}, \ \ 
\begin{math} 2^{5} \times 3^{3} \times 5^{2} \end{math}, \ \ 
\begin{math} 2^{7} \times 3^{4} \times 5^{2} \end{math}, 
\begin{math} 2^{6} \times 3^{3} \times 5^{2} \times 7^{} \end{math}, \ \ 
\begin{math} 2^{7} \times 3^{4} \times 5^{} \times 7^{} \end{math}, \ \ \begin{math} 2^{8} \times 3^{2} \times 5^{2} \times 7^{} \end{math}, \ \ \begin{math} 2^{5} \times 3^{4} \times 5^{2} \times 7^{} \end{math}, \ \ \begin{math} 2^{8} \times 3^{4} \times 5^{2} \end{math}, \ \ \begin{math} 2^{7} \times 3^{3} \times 5^{2} \times 7^{} \end{math}, \ \ \begin{math} 2^{6} \times 3^{4} \times 5^{2} \times 7^{} \end{math}$\}$.\\

One can check directly that for every integer $1\leq i \leq 3$, $S(A_i)= \dfrac{\sigma(n!)}{4}$. Thus, by Proposition \ref{3layered}, the number $n!=3628800=2^{8}\times3^{4}\times5^{2}\times7$ is a $4$-layered number.

 \textbf{Case 3}: \ $n\geq11$. Let $n\geq 11$ and $2^{\alpha_1}p^{\alpha_2}_2p^{\alpha_3}_3\cdots p^{\alpha_k}_k$ be a prime factorization of $n!$. Now let $t=(p_k-1)!$. By the definition of $n$, $(p_k-1)\geq 10$. Thus, there exists an odd almost practical number $m'$ such that $t=2^{\alpha}m'$ (see the proof of Corollary \ref{factorial 3} for  Case 5). Also, we know that $m'=p^{\beta_2}_2p^{\beta_3}_3\cdots p^{\beta_{k-1}}_{k-1}$, where $1\leq \beta_i\leq \alpha_i$ for every integer $1\leq i \leq k-1$. Therefore, by Proposition \ref{almost divisors}, the number $m=p^{\alpha_2}_2\cdots p^{\alpha_{k-1}}_{k-1}$ is an odd almost practical number. Also, by Theorem \ref{f3}, $\alpha_{k-1}=1$; this concludes that $2\mid\sigma(m')$. Then, by Proposition \ref{17}, $m$ is an odd Zumkeller number. Moreover, by Theorem \ref{f3}, $p_k<\sigma(2^{\alpha_1})$. Therefore, by Corollary \ref{2p}, the number $2^{\alpha_1}p_k$ is a Zumkeller number. Thus, by Proposition \ref{kzum}, $n!=2^{\alpha}p_km$ is a $4$-layered number.
  \end{proof}
  \begin{remark}\label{41f}
In Case 3 above, we proved that for every integer $n\geq11$ there exist Zumkeller numbers $m_1$ and $m_2$ such that $\gcd(m_1, m_2)=1$ and $n!=m_1m_2$. Without loss of generality, we let $m_2$ be odd. By Theorem 4.12 of \cite{Mahanta}, the number $2^{\ell}m_1$ is a Zumkeller number for every positive integer $\ell$. Thus, by Proposition \ref{nw}, the number $n!$ is a good $4$-layered number for $n\geq11$. Also, by Proposition \ref{pqr}, the number $9!$ is a good $4$-layered number. In addition, one can check directly that the numbers $10!$ is good $4$-layered; we omit the details here. Therefore, we can say that the number $n!$ is $4$-layered if and only if $n!$ is a good $4$-layered number.
\end{remark}
\begin{remark}\label{42f}
  By Theorem \ref{f3} and Corollary \ref{2p}, $2^{\floor*{\frac{n}{2}}}p_k$ is also a Zumkeller number. Thus, using the same arguments of Case 3 of the above proof, we get $\dfrac{n!}{2^t}$ is a 4-layered number for $n\geq 11$, where $1\leq t\leq \sum_{i=2}^{\infty}\floor*{\dfrac{n}{2^i}}$. Therefore, by the previous remark, if $r=\sum_{i=2}^{\infty}\floor*{\dfrac{n}{2^i}}$, then the number $\dfrac{n!}{2^r}$ is a good $4$-layered number for $n\geq11$.  
  \end{remark}
      \begin{remark}\label{faczum}
By Bertrand postulate and Proposition \ref{prop5}, every primorial number \footnote{For the $n$th prime number $p_n$, the primorial $p_n\#$ is defined as the product of the first $n$ primes.} $n\neq2$ is a practical number. Also, by Proposition \ref{prop5}, if $n$ is a practical number and $d\mid n$, then $nd$ is also a practical number. Then, By Proposition \ref{17}, if $n=p^{\alpha_1}_1p^{\alpha_2}_2\cdots p^{\alpha_{\ell}}_{\ell}$, where $p_i$ be the $i$th prime, and $\alpha_1, \alpha_2,\dots, \alpha_{\ell}$ are positive integers for every positive $1\leq i \leq \ell$, then the number $n$ is a Zumkeller number if and only if $2\mid \sigma(n)$.  In addition, by  Bertrand's postulate and the same method that we used to prove Theorem \ref{f3}, we can prove that for integers $n\geq3$, there exists a prime number $q_n$ such that $q_n\mid\mid n!$. Thus, $2\mid\sigma(n!)$ for integers $n\geq3$. Therefore, by Proposition \ref{17}, $n!$ is a Zumkeller number for integers $n\geq3$.
\end{remark}
 \begin{remark}\label{43f}
Let $n\leq 10^{10}$ be a practical number. One can check directly that if is good $k$-abundant, then $n$ is $k$-layered.
  \end{remark}
  Proposition \ref{prop5}, Theorems \ref{factorial 3} and \ref{factorial 4}, Remarks \ref{38f},  \ref{39f}, \ref{41f}, \ref{42f}, and \ref{43f} leads us to the following conjecture.
\begin{conjecture}\label{factorial conjecture}
Let $n$ be a practical number. If $n$ is good $k$-abundant, then $n$ is $k$-layered.
  \end{conjecture}
Now we recall a well-known proposition which is a consequence of prime number for arithmetic progressions.
 \begin{proposition}
Let $m$ be a positive integer. There exists a positive integer $s_m$ such that for integers $n\geq s_m$, there exists a prime $p_{m,n}$ between $n$ and $2n$, where $p_{m, n} \equiv m-1 \pmod m$.
 \end{proposition}
  By using the previous proposition and the same method that we used in the proof of Theorem \ref{f3}, we can conclude the following theorem.
  \begin{theorem}\label{mohem}
 Let $m$ be a positive integer. For every positive integer $n$ that is large enough, then there exists a prime $p_{m, n}$ such that $p_{m, n} \equiv m-1 \pmod m$ and $p_{m, n} \mid\mid n!$. 
  \end{theorem}
    Now we state a well-known proposition which recalls some elementary properties of the abundancy function $I$.
 \begin{proposition}\label{I(n)}
  The abundancy function $I$ possesses following properties.
  
  $\mathrm{(i)}\ $ Let $p$ and $q$ be primes such that $p<q$. Let $\alpha$ and $\beta$ be positive integers such that $\beta \leq \alpha$. Then, $I(q^{\beta})<I(p^{\alpha}).$
  
 $\mathrm{(ii)}\ $ Let $p^{\alpha_1}_1 p^{\alpha_2}_2\cdots p^{\alpha_k}_k$ be the prime factorization of a positive integer $n$. Then
  $$I(n)=\prod_{i=1}^{k}I(p^{\alpha_i}_i)<\prod_{i=1}^{k} \dfrac{p_{i}}{p_{i}-1}.$$
  \end{proposition}
Thus, by the Theorem \ref{mohem} and some elementary properties of abundancy function, for every positive integer $k$, there exists a positive integer $s_k$ such that $k\mid \sigma(n!)$ and $I(n!)\geq k$ for integers $n\geq s_k$. Then we have the following corollary.
  \begin{corollary} 
  Let $k$ be a positive integer. If Conjecture \ref{factorial conjecture} holds, then there exists a positive integer $s_k$ such that for integers $n\geq s$, the number $n!$ is a good $k$-layered number.
  \end{corollary}
  Let $\alpha$ and $1\leq \beta \leq 3$ be positive integers. We close this section by the following theorem which find all $4$-layered numbers of the form $n=p^{\alpha}_1p^{\beta}_2p_3p_4$, where $p_i$ is prime for integers $1\leq i\leq 4$.
 \begin{theorem}\label{pqr}
 Let $\alpha$ and $1\leq \beta \leq 3$ be positive integers. Also, let $n=p^{\alpha}_1p^{\beta}_2p_3p_4$, where $p_i$ is prime for integers $1\leq i\leq 4$. The number $n$ is 4-layered if and only if $p_1=2, p_2=3, p_3=5, p_4=7, \alpha \geq 5,$ and $\beta=3$
 \end{theorem}
 \begin{proof}
 By some properties of the abundancy function, we know that $I(2^{\alpha}3^{\beta}5^{2})<\dfrac{2}{2-1}\times\dfrac{3}{3-1}\times\dfrac{1+5+5^2}{5^2}<4\leq I(p^{\alpha}_1p^{\beta}_2p_3p_4)$. Thus, by Proposition \ref{I(n)}, for integers $1\leq i\neq j\leq 4$, $p_i\neq p_j$. Now we prove that $2^5\mid n$ by contradiction. We let $2^5\nmid n$. Therefore, $I(n)\leq I(2^5\times 3^{\alpha}\times 5\times 7)<4$, a contradiction. This concludes that $p_1=2$ and $\alpha\geq 4$. Also, $3^{3}\mid n$, for otherwise $I(n)\leq I(2^{\alpha}\times 3^{2}\times 5\times 7)<4$, a contradiction. By the same method, one can prove that $p_3=5$, $p_4=7$.
 Thus, for completing the proof, it is sufficient to proof that the number $n=2^{5}\times3^{3}\times5\times7$ is a good $4$-layered number. Let 
 \begin{longlist}
$A_{1} = \{$ \begin{math} 2^{5} \times 3^{3} \times 5^{} \times 7^{} \end{math}$\},$
\end{longlist}

\begin{longlist}
$A_{2} = \{$ \begin{math} 2^{} \times 3^{} \end{math}, \ \ 
\begin{math} 2^{4} \times 3^{2} \end{math}, \ \ 
\begin{math} 2^{5} \times 3^{} \times 5^{} \end{math}, \ \ 
\begin{math} 2^{} \times 3^{3} \times 5^{} \times 7^{} \end{math}, \ \ 
\begin{math} 2^{4} \times 3^{2} \times 5^{} \times 7^{} \end{math}, \ \ 
\begin{math} 2^{3} \times 3^{3} \times 5^{} \times 7^{} \end{math}, \ \ 
\begin{math} 2^{4} \times 3^{3} \times 5^{} \times 7^{} \end{math}$\},$
\end{longlist}

\begin{longlist}
$A_{3} = \{$ \begin{math} 3^{} \end{math}, \ \ 
\begin{math} 2^{4} \times 3^{} \end{math}, \ \ 
\begin{math} 3^{3} \times 5^{} \times 7^{} \end{math}, \ \ 
\begin{math} 2^{4} \times 3^{3} \times 5^{} \end{math}, \ \ 
\begin{math} 2^{3} \times 3^{2} \times 5^{} \times 7^{} \end{math}, \ \ 
\begin{math} 2^{4} \times 3^{3} \times 7^{} \end{math}, \ \ 
\begin{math} 2^{5} \times 3^{} \times 5^{} \times 7^{} \end{math}, \ \ 
\begin{math} 2^{2} \times 3^{3} \times 5^{} \times 7^{} \end{math}, \ \ 
\begin{math} 2^{5} \times 3^{3} \times 5^{} \end{math}, \ \ 
\begin{math} 2^{5} \times 3^{2} \times 5^{} \times 7^{} \end{math}$\}$,
\end{longlist}
\noindent and $A_4= D_n\setminus (A_1 \cup A_2 \cup A_3)$. One can check directly that $\{A_1, A_2, A_3, A_4\}$ is a $4$-layered partition for $n$. Let $n=2^{5}m$. Let 
$A'_1 = \{$ \begin{math} 3^{3} \times 5^{} \times 7^{} \end{math}, \ \ 
\begin{math} 3^{} \times 5^{} \end{math}$\}$ and $A'_2=D_m\setminus A'_1$. One can check directly that $\{A'_1, A'_2\}$ is a Zumkeller partition for $m$. Also, if we define $A''_1=\{2^{5}d : d\in A'_1\}$ and $A''_2=\{2^{5}d : d\in A'_2\}$, then one can check directly that $A''_1\subset (A_1\cup A_2)$ and $A''_2\subset (A_3\cup A_4)$. Thus, by Proposition \ref{is 4layered}, $n$ is a good $4$-layered number.
 \end{proof}
\section{$k$-multiperfect numbers and $k$-layered numbers}\label{multiperfect number section}
 We recall that a number $n$ is said to be $k$-multiperfect if $\sigma(n)=kn$ for some integer $k\geq 2$. It is clear that if $n$ is $2$-multiperfect, then $n$ is said to be perfect. Since every perfect number is Zumkeller, this suggests us the following open question.
\begin{open question}\label{open1}
Is every $k$-multiperfect number  $k$-layered?
 \end{open question}
 In the following, we will address this open question in some special cases. Before that, we recall that a number $n$ is said to be semiperfect if $n$ is equal to the sum of all or some of its proper divisors. Also, the abundant number $n$ which is not semiperfect is called weird. The existence of odd weird numbers is still an open question. 
The number $70$ is the smallest weird number. If $n$ is weird and $p$ is a prime number such that $\sigma(n)<p$, then $np$ is a weird number (for more details, see \cite{semiperfect2}). This shows that there exist infinity many weird numbers. Nevertheless, it is not yet known whether there exists infinitely many primitive weird numbers \footnote{A weird number $n$ that is not multiple of other weird number is said to be primitive weird number.} or not. Recently, a wide range of number theorists have tried to introduce some optimal algorithms to find new primitive weird numbers. (For instance, see \cite{melfi2, melfi1}.)
Fang and Beckert proved, using parallel tree search, that there are no odd weird numbers up to $10^{21}$(\cite{fang}). This leads us to the following conjecture.
 \begin{conjecture}\label{semiperfect}
 Every odd abundant number is semiperfect.
 \end{conjecture}
Now we recall two well-known propositions about semiperfect numbers.
 \begin{proposition}\label{semi}
 Let $m$ be a positive integer, and let $p$ be a odd prime number such that $p \leq 2^{m+1}-1$. Then, the number $n=2^{m}p$ is a semiperfect number.
 \end{proposition}
\begin{proposition}\label{semi2}
 If a positive integer $m$ is divisible by a semiperfect number $n$, then $m$ is semiperfect.  
\end{proposition}
 \begin{remark}
 Let $\alpha$ and $\beta$ be positive integers. Also, let $p_1$ and $p_2$ be primes. By Remark \ref{Pankaj}, the number $n=p_1^{\alpha_1}p_2^{\alpha_2}$ is an abundant number if and only if $p_1=2$ and $3\leq p_2\leq 2^{\alpha+1}-1$. Also, by Propositions \ref{semi} and  \ref{semi2}, all numbers of this form are semiperfect. This concludes that all abundant number which possess  exactly two distinct prime factors are semiperfect. 
 \end{remark}
It is believed that all $k$-multiperfect number of abundancy $3, 4, 5, 6$ and $7$ are known. Reference \cite{list} contains the list of all known $k$-multiperfect numbers. According to this list, one can directly obtain the following proposition.

\begin{proposition}\label{3perfect}

$\mathrm{(i)}\ $ All the known $k$-multiperfect numbers are practical.

 $\mathrm{(ii)}\ $ All the $6$ known $3$-multiperfect numbers are divisible by $3$ or $5$. Also, all $3$-multiperfect numbers are divisible by $2^{3}$.
 
 $\mathrm{(iii)}\ $ All the $36$ known $4$-multiperfect numbers are divisible by $2^{2}$. In addition, if $n=2^{\alpha}t$ is one the known $4$-multiperfect, then $\dfrac{n}{2}$ satisfies the conditions of Proposition \ref{prop5}; this concludes that $\dfrac{n}{2}$ is a practical number.
 
 $\mathrm{(iv)}\ $  Let $n$ be one of the known $5$-multiperfect numbers. Then, $2^{7}\mid n$ and $5\mid n$ or $7\mid n$.
\end{proposition}
As a consequence of part $\mathrm{(i)}\ $ of the above proposition, we have the following corollary.
\begin{corollary}
If Conjecture \ref{factorial conjecture} holds, then every known $k$-multiperfect number is $k$-layered.
\end{corollary}
Thus, we have the following theorem.
\begin{theorem}\label{every3perfect}
Every known $3$-multiperfect number is $3$-layered.
\end{theorem}
\begin{proof}
Let $n$ be one of the known $3$-multiperfect numbers. By Propositions \ref{semi} and \ref{semi2}, and part $\mathrm{(ii)}\ $ of Proposition \ref{3perfect}, the number $n$ is semiperfect. This concludes that there exists  $A_1\subset D_n\setminus\{n\}$ such that $S(A_1)=n=\dfrac{\sigma(n)}{3}$. Now we define $A_2=\{n\}$ and $A_3=D_n\setminus(A_1\cup A_2)$. Thus, by the definition of $n$, $\{A_1, A_2, A_3\}$ is a $3$-layered partition for $n$.
\end{proof}
Now we state an essential and simple lemma. 
\begin{lemma}\label{4mul}
Let $n=2^{\alpha}t$ is an even positive integer with $t$ odd and $I(n)\geq 4 $. Then, $I(t)> 2$.
\end{lemma}
\begin{proof}
We know that $I(2^{\alpha})=\dfrac{2^{\alpha+1}-1}{2^{\alpha}}<2$. This concludes that $I(t)>2$.
\end{proof}
Now by assuming Conjecture \ref{semiperfect}, we state a theorem similar to Theorem \ref{every3perfect} about all the known $4$-multiperfect numbers.
 \begin{theorem}\label{4con}
 If Conjecture \ref{semiperfect} holds, then every known $4$-multiperfect numbers are $4$-layered.
 \end{theorem}
 \begin{proof}
 Let $n=2^{\alpha}m$ be a known $4$-multiperfect number with $\gcd(2^{\alpha}, m)=1$.
 By part $\mathrm{(iii)}\ $ of Proposition \ref{3perfect}, the number $t=\dfrac{n}{2}$ is a practical number. Also, by Lemma \ref{4mul}, $\sigma(t)>n$. Thus, by definition of practical numbers, there exists $A_1\subset D_t\setminus\{t\}$ such that $S(A_1)=n$=$\dfrac{\sigma(n)}{4}$. Moreover, by assuming  Conjecture \ref{semiperfect}, there exists $A\subset D_m\setminus\{m\}$ such that $S(A)=m$. Now we define $A_2=\{2^{\alpha}d: d\in A\}$. It is clear that $S(A_2)=n$=$\dfrac{\sigma(n)}{4}$. Furthermore, we define $A_3=\{n\}$. By definition of $A_1, A_2,$ and $A_3$ for every integers $1\leq i \neq j \leq 3 $,
  $A_i\cap A_j=\emptyset$. Then, by Proposition \ref{3layered}, $n$ is $4$-layered.
 \end{proof}
 Now we proof a theorem similar to Theorem \ref{4con} about $5$-multiperfect numbers.
 \begin{theorem}\label{5mul}
By assuming Conjecture \ref{semiperfect}, every known $5$-multiperfect numbers are $5$-layered.
 \end{theorem}
 \begin{proof}
 Let $n$ be one of the known $5$-multiperfect numbers. By part  $\mathrm{(iv)}\ $ of Proposition \ref{3perfect}, we know there exists a prime number $p=5$ or $p=7$ where $p\mid n$. Now let $\alpha$, $\beta$, $t$, and $m$ be positive integers such that $\nu_2(n)=\alpha$, $\nu_p(n)=\beta$, and $n=2^{\alpha}p^{\beta}t=mt$. we have
 $I(m)<\dfrac{2}{2-1}\times\dfrac{5}{5-1}=2.5$.

We know $I(n)=I(m)I(t)=5$; this concludes that $I(t)>2$. By Proposition \ref{semi}, there exists $A\subset D_m\setminus\{m\}$ where $S(A)=m$. Also, by assuming Conjecture \ref{semiperfect}, there exists $A'\subset D_t\setminus\{t\}$ where $S(A')=t$. Now we define
$A_1=\{n\}, A_2=\{dt : d\in A\}, A_3=\{dm : d\in A' \},$ and $A_4=\{d_1d_2 : d_1\in A \land d_2\in A'\}$.
One can check directly that for every integer $1\leq i\leq 4$, $S(A_i)=\dfrac{\sigma(n)}{5}=n$. Thus, by Proposition \ref{3layered}, $n$ is $5$-layered.
 \end{proof}
 \begin{remark}
For every known 4-multiperfect number $n$, one can check that the number $m$, which is defined in the proof of Theorem \ref{4mul}, is semiperfect. In addition, for every known 5-multiperfect number $n$, one can check that the number $t$, which is defined in the proof of Theorem \ref{5mul}, is semiperfect. This concludes that the all known 4-multiperfect numbers and the all known 5-multiperfect numbers are 4-layered and 5-layered, respectively. We do not pursue it here.
 \end{remark}
 Now we state an essential proposition.
 \begin{proposition}\label{f4}
  Let $n$ and $m$ be positive integers with $gcd(n, m)$=1. If $n$ is $k$-layered and $\dfrac{\sigma(m)}{k+1}$ is a sum of a subset of the set of all positive divisors of $m$, then $nm$ is $(k+1)$-layered.
 \end{proposition}
 \begin{proof}
Let $\{A_1, A_2, \dots, A_k\}$ be a $k$-layered partition for $n$. Also, let  $A$ be a subset of the set of all positive divisors of $m$ such that $S(A)=\dfrac{\sigma(m)}{k+1}$. We define $A'=D_m\setminus A$ One can check directly that for every integer $1\leq i\leq k$, $S(A_iA')=\dfrac{\sigma(nm)}{k+1}$. Thus, by Proposition \ref{3layered}, $nm$ is $(k+1)$-layered. 
 \end{proof}
As a consequence of Proposition \ref{f4}, we have the following corollary.
\begin{corollary}\label{f5}
Let $p$ and $n$ be a prime number and a $p$-layered number, respectively. If $\gcd(n, p)=1$, then $pn$ is a $(p+1)$-layered. 
\end{corollary} 
\begin{remark}\label{f6}
Let $p$ and $n$ be a prime number and a $p$-multiperfect number, respectively. One can check directly that if $\gcd(n, p)=1$, then $pn$ is a $(p+1)$-multiperfect number. \end{remark}
By Corollary \ref{f5} and Remark \ref{f6}, we have the following theorem.
\begin{theorem}\label{zaman}
Let $p$ be a prime number and $n$ be a positive integer such that $\gcd(n, p)=1$. If $n$ is $p$-multiperfect and $p$-layered, then $pn$ is $(p+1)$-multiperfect and $(p+1)$-layered.
\end{theorem}
The two following examples state a application of previous theorem to find a new family of $k$-multiperfect numbers which are $k$-layered.
\begin{example}
 One can check that the number $n=2^{14}\times5\times7\times19\times31\times151$ is a $3$-multiperfect number. By Theorem \ref{every3perfect}, $n$ is $3$-layered. Thus, by Theorem \ref{zaman}, $3n$ is $4$-multiperfect and $4$-layered
\end{example}
\begin{example}
 One can check that the number $n=2^{29}\times3^{10}\times7^{5}\times11^{2}\times13\times19^{3}\times23\times31\times43\times83\times107\times151\times181\times331\times3851$ is a $5$-multiperfect number. By Remark \ref{f6}, $5n$ is $6$-multiperfect. If conjecture \ref{semiperfect} holds, then by Theorem \ref{4con}, $n$ is $5$-layered. Also, if $n$ is $5$-layered, then by Corollary \ref{f5}, $5n$ is $6$-layered.
\end{example}
Now we state a proposition by which we can find more family $k$-multiperfect numbers which are $k$-layered.
\begin{proposition}\label{2013}
Let $n$ and $m$ be a $k$-multiperfect number and a $(k+1)$-multiperfect number, respectively, such that $n\mid m$. If $n$ is a $k$-layered number, then $m$ is $(k+1)$-layered number.
\end{proposition}
\begin{proof}
Let $\{A_1, A_2, \dots, A_k\}$ be a $k$-layered partition for the $k$-multiperfect number $n$. Now let  $dn=m$ . For every integer $1\leq i \leq k$, we define $B_i=\{ad: a\in A_i\}$. Then, for every integer $1\leq i \leq k$, $\sigma(B_i)=m=\dfrac{\sigma(m)}{k+1}$. Thus, by Proposition \ref{3layered}, $m$ is $(k+1)$-layered.
\end{proof}
\begin{remark}
One can check directly that exactly half of known $4$-multiperfect numbers are divisible by at least  a $3$-multiperfect number. Then, by Propositions \ref{every3perfect} and \ref{2013}, at least half of known $4$-multiperfect are $4$-layered.
\end{remark}
We close this section by an example that implies an application of Proposition \ref{2013} to find another new family of $k$-multiperfect number which are $k$-layered.
\begin{example}
 Let $t_1=2\times3, t_2=2^3\times3\times5, t_3=2^5\times 3^3\times 5 \times 7, t_4=2^{11}\times 3^3\times 5^2\times 7^2\times 13\times 19\times 31, t_5=2^{19}\times 3^5\times 5^2\times 7^2\times 11\times 13^2\times 19^2\times 31^2\times 37\times 41\times 61\times 127, t_6=2^{39}\times 3^{11}\times 5^7\times 7^3\times 11\times 13^2\times 17\times 19^2\times 29\times 31^2\times 37\times 41\times 61\times 73\times 79\times 83\times 127\times 157\times 313\times 331\times 2203\times 30841\times 61681$. One can check that for integers $1\leq i\leq 6$, the number $t_i$ is a $(i+1)$-perfect number. Also, it is easy to see that for every integer $1\leq i\leq 5$, $t_i\mid t_{i+1}$. We know that $6$ is Zumkeller. Thus, by Proposition \ref{2013}, for every integer $1\leq i\leq 6$, $t_i$ is a $(i+1)$-layered number. 
\end{example}
\section{Near-perfect numbers and Zumkeller numbers}\label{near perfect number section}
In this short section, we do not apply any advanced method. In fact, the importance of Remark \ref{nearremark} motives the author to state this section; by this remark, we can find a new family of Zumkeller numbers.

We recall that a number $n$ is said to be a near-perfect number if $n$ is the sum of all of its proper divisors, except for one of them, which we term the redundant divisor \cite{Pollack}. By definition of near-perfect numbers, we have the following well-known proposition.
\begin{proposition}\label{clear near}
A number $n$ is a near-perfect number with redundant divisor $d$ if and only if $d$ is a proper divisor of $n$, and $\sigma(n)=2n+d$
\end{proposition}
\begin{remark}\label{odd near-perfect}
By Proposition \ref{clear near}, if a number $n$ is an odd near-perfect number, then $\sigma(n)$ is odd. Thus, by Proposition \ref{clear}, $n$ fails to be Zumkeller.
\end{remark}
Now we state a proposition similar to Proposition \ref{17}, about near-perfect numbers.
\begin{proposition}\label{near}
Let $n$ be a near-perfect number. Then, $n$ is Zumkeller if and only if $\sigma(n)$ is even.
\end{proposition}
\begin{proof}
By Proposition \ref{clear}, if $n$ is a Zumkeller number, then $\sigma(n)$ is even. Now let $\sigma(n)$ is even. By the definition of $n$ and Proposition \ref{clear near}, $n$ possesses a divisor $d$ such that $2d$ is a proper divisor of $n$ and $\sigma(n)=2n+2d$. We define $A=\{n, d\}$. One can check directly that $\{A, D_n\setminus A \}$ is a Zumkeller partition for $n$.
\end{proof}
The following proposition was stated by Ren and Chen \cite{Ren}; this proposition classifies all near-perfect numbers with two distinct prime factors.
\begin{proposition}\label{class}
Let $n$ be a positive integer with two prime factors. The number $n$ is near-perfect if and only if one of the following holds:

$\mathrm{(i)}\ n=2^{t-1}(2^{t}-2^{k}-1)$, where $2^{t}-2^{k}-1$ is prime;

$\mathrm{(ii)}\ n=2^{2p-1}(2^p-1)$, where $2^p-1$ is a Mersenne prime;

$\mathrm{(iii)}\ n=2^{p-1}(2^{p}-1)^{2}$, where $2^p-1$ is a Mersenne prime;

$\mathrm{(iv)}\ n=40$.

\end{proposition}

As a direct consequence of Propositions \ref{near} and \ref{class}, we have the following corollary.
\begin{corollary}
Let $n$ be a near-perfect number with two distinct prime factors. The number $n$ is non-Zumkeller if and only if there exists a Mersenne prime $p$ such that $n=2^{p-1}(2^{p}-1)^{2}$.
\end{corollary}
Now we recall a proposition of \cite{Li}. This proposition classifies all near-perfect numbers  of the form $2^{\alpha}p_1p_2$, where $p_1$ and $p_2$ are odd primes with $p_1<p_2$.
\begin{proposition}\label{li and liao}
Let $\alpha$ be a positive integer. Let $p_1$ and $p_2$ be odd primes with $p_1<p_2$. Then, the number $n=2^{\alpha}p_1p_2$ is near-perfect if and only if $n$ satisfies one the following conditions.

$\mathrm{(i)}\ p_1=\dfrac{2^{\alpha+1}-1+k}{2^{\beta}+1-k}$, where $k=\dfrac{2^{\alpha+1}-1}{p_2}$ and $1\leq \beta \leq \alpha-1$. 

$\mathrm{(ii)}\ $ $p_1=2^{\alpha+1}-1+\dfrac{2^{\alpha}-2^{\beta-1}}{k}$, where $k$ is determined by equation $$p_2=(2^{\alpha+1}-1)(2k+1)-2^{\beta}, 1 \leq \beta \leq \alpha.$$

$\mathrm{(iii)}\ p_2=2^{\alpha+1}-1 + \dfrac{2^{2\alpha+1}-2^{\alpha}-2^{\beta -1}}{k}$, where $k=\dfrac{p_1-(2^{\alpha+1}-1)}{2}$ and $1\leq \beta \leq \alpha$.
\end{proposition}
\begin{remark}\label{nearremark}
Let $n=2^{\alpha}p_1p_2$ be a number satisfying case $\mathrm{(i)}\ $ of Proposition \ref{li and liao}. One can check directly that $p_1<2^{\alpha+1}-1$. This concludes that $2^{\alpha}p_1$ is a practical number. Thus, without knowing Propositions \ref{near} and \ref{li and liao}, by Propositions \ref{prop5} and \ref{17}, we know that all the numbers satisfying case $\mathrm{(i)}\ $ of Proposition \ref{li and liao} are Zumekeller, but by case $\mathrm{(ii)}\ $and case $\mathrm{(iii)}\ $ of Proposition \ref{li and liao}, we can construct a new subset of Zumkeller numbers. 
\end{remark}
In \cite{Li}, the authors checked that there exists exactly $8$ near-perfect numbers satisfying case $\mathrm{(i)}\ $of Proposition \ref{li and liao} for $\alpha<1000$, $\beta=\dfrac{\alpha+5}{2}$, and $k=2^{\beta-2}+1$; the number $2^9 \times 11 \times 31$ is the smallest of such numbers. Also, for case $\mathrm{(ii)}\ $ of Proposition \ref{li and liao}, They checked that there exist exactly $289$ near-perfect number; the number $2^{3}\times17\times101$ is the smallest such number. Furthermore, for case $\mathrm{(iii)}\ $ of Proposition \ref{li and liao}, they checked that there exist exactly $248$ near-perfect number for $\alpha<100$ and $k\leq 100$; the number $2^{3}\times17\times131$ is the smallest such number.

In \cite{clark}, the authors raised the following conjecture.
\begin{conjecture}\label{zumkeler conjecture}
If $n$ is an even and Zumkeller number, then $n$ is half-Zumkeller.
\end{conjecture}
In \cite{zumkeller}, the authors verified that the conjecture is true in some cases. They proved that if $n$ is an even Zumkeller number such that $\sigma(n)<3n$, then $n$ is a half-Zumkeller. Thus, by Remark \ref{odd near-perfect}, we have the following proposition.
\begin{proposition}
Let $n$ be a near-perfect number. If $n$ is a Zumkeller number, then $n$ is half-Zumkeller.
\end{proposition}

\section{The difference of two consecutive $k$-layered numbers}\label{difference section}
 In this section, we let $ k, m, a, b, s, r$, and $z$ be positive integers. We also let $p_1, p_2, \cdots, p_m$ be distinct primes such that $p_1<p_2<\dots<p_m$. We start this section with two new definitions.
\begin{definition}
 The ascending chain $\ell=(a+t)^{z }_{t=0}$ is said to be a $(p_1, p_2, \dots, p_m)$-$\gcd$ chain if for every integer $0\leq j\leq z$, $\gcd(a+j, p_1p_2\cdots p_m)\neq1$.
\end{definition}
\begin{definition}
The chain $\ell=(a+t)^{z}_{t=0}$ is said to be a $(p_1, p_2, \dots,  p_m)$-$\gcd$-max chain if $\ell$ satisfies the following conditions:

$\mathrm{(i)}\  \ell$ is a $(p_1, p_2, \dots, p_m)$-$\gcd$ chain;

$\mathrm{(ii)}\ $ If $\ell'=(a+t)^{s}_{t=0}$ is a 
$(p_1, p_2, \dots, p_m)$-$\gcd$ chain, then $z\geq s$.

We also define $L_{max}(p_1, p_2, \dots, p_m)$ as the cardinality of $\ell$.
\end{definition}

We recall a standard well-known proposition, which is a generalization of the Chinese reminder theorem.
\begin{proposition}\label{chinese}
 Let $n_1, n_2, \dots, n_k, t_1, t_2, \dots,$ and $t_k$ be positive integers. Then, the simultaneous set of congruences  $x \equiv t_1 \pmod {n_1}$, $x \equiv t_2 \pmod {n_2}$, $\dots$, and  $x \equiv t_k \pmod {n_k}$ has a solution if and only if for every integers $1\leq i\neq j\leq k$, $\gcd(n_i, n_j)\mid t_i-t_j$.
\end{proposition}
The following proposition finds a lower bound for $L_{max}(p_1, p_2, \dots, p_m)$.
\begin{proposition}\label{lower}
 $L_{max}(p_1, p_2, \dots, p_m)\geq m$.
\end{proposition}
\begin{proof}
By Proposition \ref{chinese}, we can find a positive integer $a$ such that $a \equiv 0$ (mod $p_1$), $a \equiv -1$ (mod $p_2$), $\dots$, and $a \equiv -(m-1)$ (mod $p_m$). Then, the chain $\ell=(a+t)^{m-1}_{t=0}$ is a $(p_1, p_2, \dots, p_m)$-$\gcd$ chain.
\end{proof}

We know that $L_{max}$ can be known as a function defined on the set of all positive square-free integers. The following proposition state another property of this function.
\begin{proposition}
 $L_{max}(p_1, p_2, \dots, p_m)<L_{max}(p_1, p_2, \dots, p_m, p_{m+1})$.
\end{proposition}
\begin{proof}
Let $\ell=(a+t)^{s}_{t=0}$ be a $(p_1, p_2, \dots, p_m)$-$\gcd$ chain. This concludes that for every integer $0\leq t \leq s$, there exists positive integer $1\leq i_t \leq m$ such that $p_{i_t}\mid a+t$. Thus, by Proposition \ref{chinese}, if $0\leq t \neq t' \leq s$, then $\gcd(p_{i_t}, p_{i_{t'}})\mid t-t' $. Therefore, by the definition of $p_{m+1}$ and Proposition \ref{chinese}, there exists a positive integer $b$  such that  $b \equiv 0$ (mod $p_{i_0}$), $b \equiv -1$ (mod $p_{i_1}$), $\dots$, $b \equiv -s$ (mod $p_{i_s}$), and $b \equiv -(s+1)$ (mod $p_{m+1}$). Thus, $\ell'=(b+t)^{s+1}_{t=0}$ is a $(p_1, p_2, \dots, p_m, p_{m+1})$-$\gcd$ chain. This completes the proof.
\end{proof}

In remainder of this note, we define $\phi$ to be the Euler's totient function. Let $n$ be a positive integer with prime factorization $p^{\alpha_1}_1p^{\alpha_2}_2\cdots p^{\alpha_m}_m$. We recall that $\phi(n)=n\prod_{i=1}^{m}(1-\dfrac{1}{p_i})$. In Proposition \ref{lower}, we found a lower bound for $L_{max}(p_1, p_2, \dots, p_m)$. Now we find an upper bound for $L_{max}(p_1, p_2, \dots, p_m)$.

\begin{proposition}\label{upper for L}
$L_{max}(p_1,p_2,\dots, p_m)\leq p_1p_2\cdots p_m-(p_1-1)(p_2-1)\cdots(p_m-1).$
\end{proposition}
\begin{proof}
Let $t=p_1p_2\cdots p_m$. Let $r$ be a non-negative integer such that $0\leq r<t$. We know that for every positive integer $j$, $\gcd(jt+r, t)=1$ if and only if $\gcd(r, t)=1$. Therefore, between $0$ and $t$ there exist exactly $\phi(t)$ distinct integers $r_1, r_2, \dots,$ and $r_{\phi(t)}$ such that $1 = r_1 < r_2 < \dots < r_{\phi(t)}< t$ and  $\gcd(jt+r_i, t)=1$ for every integers $1\leq i\leq \phi(t)$ and $ j\geq 0 $. This concludes that between $t$ and $2t$ there exist exactly $\phi(t)$ positive integers $n_1< n_2< \dots< n_{\phi(t)}$ such that  $\gcd(n_i, t)=1$  for every integer $1\leq i \leq \phi(t)$. Thus, $L_{max}(p_1,p_2,\dots,p_m)\leq t-\phi(t).$ This completes the proof. 
\end{proof}

The following proposition calculates  the function $L_{max}$  at some value.
\begin{proposition}\label{L at some value}
$\mathrm{(i)}\ L_{max}(3, 5, 7, 11)=6$.

$\mathrm{(ii)}\ L_{max}(3, 5, 7, 11, 13, 17)=12$.
\end{proposition}
\begin{proof}
$\mathrm{(i)}\ $ By the Chinese remainder theorem, there exists a positive integer $a$ such that $a \equiv 0$ (mod $5$), $a \equiv -1$ (mod $3$), $a \equiv -2$ (mod $7$), and $a \equiv -3$ (mod $11$). So, $a\equiv -4$ (mod $3$) and $a\equiv -5$ (mod $5$). Therefore, $\ell=(a+t)^{5}_{t=0}$ is a $(3, 5, 7, 11)$ - $\gcd$ chain. Now it is sufficient to prove that for every positive integer $a'$, there exists an integer $0\leq i \leq 6$ such that $\gcd(a'+i, 3\times5\times7\times11)=1$. Let $\ell'=(a'+t)^{6}_{t=0}$ be a chain of positive integers, and let $D=\{3, 5, 7, 11\}$. For every integer $d\in D$, We define $A_d$ as the set of all numbers $b$ of chain $\ell'$ such that $d\mid b$. Then, $|A_d|\leq \ceil*{\dfrac{7}{d}} $. We investigate the two cases for $a' $:

\textbf{Case 1}: \ $3\mid a'$. We know that $|A_5|\leq2$. Also, one can check directly that in this case, if $|A_{5}|=2$, then $|A_{5}\cap A_{3}|=1$. This concludes that $|A_3\cup A_5|\leq 4$. Thus, $|A_3 \cup A_5 \cup A_7 \cup A_{11} |\leq 6$.

\textbf{Case 2}: \ $3\nmid a'$. In this case, $|A_3|\leq 2$. Thus, one can check directly that  $|A_3 \cup A_5 \cup A_7 \cup A_{11} |\leq 6$.

This completes the proof.

$\mathrm{(ii)}\ $ By Chinese remainder theorem, there exists a positive integer $a$ such that $a \equiv 0$ (mod $11$), $a \equiv -1$ (mod $3$), $a \equiv -2$ (mod $7$), $a \equiv -3$ (mod $5$), $a \equiv -5$ (mod $13$), and $a \equiv -6$ (mod $17$). One can check directly that $\ell=(a+t)^{11}_{t=0}$ is a $(3, 5, 7, 11, 13, 17)$ - $\gcd$ chain.
Now it is sufficient to prove that for every positive integer $a'$, there exists an integer $0\leq i \leq 12$ such that $\gcd(a'+i, 3\times5\times7\times11\times13\times17)=1$. Let $\ell'=(a'+t)^{12}_{t=0}$ be a chain of positive integers, and let $D=\{3, 5, 7, 11, 13, 17\}$. For every integer $d\in D$, We define $A_d$ as the set of all numbers $b$ of chain $\ell'$ such that $d\mid b$. Then, $|A_d|\leq \ceil*{\dfrac{13}{d}} $. We investigate the two cases for $a' $:

\textbf{Case 1}: \ $3\mid a'$. We know that $|A_{11}|\leq 2$. Also, one can check directly that in this case, if  $|A_{11}|=2$, then $|A_{11}\cap A_{3}|=1$. Then, $|A_{3} \cup A_{11}|\leq 6$. Also, we know that $|A_{5}|\leq 3$, and if $|A_5|=3$, then $|A_{3}\cap A_5|=1$. This concludes that $|A_3\cup A_5 \cup A_{11}|\leq 8$. Thus, $|A_3 \cup A_5 \cup A_7 \cup A_{11} \cup A_{13} \cup A_{17}|\leq 12$.

\textbf{Case 2}: \ $3\nmid a'$. By the definition of $a'$, $|A_3|\leq 4$. Once again, we know that if $|A_5|=3$, then $|A_{3}\cap A_5|=1$. Therefore, $|A_3\cup A_5|\leq 6$.  Thus, $|A_3 \cup A_5 \cup A_7 \cup A_{11} \cup A_{13} \cup A_{17}|\leq 12$.

This completes the proof.
\end{proof}
We can find some more results about the function $L_{max}$. But we do not pursue this here. Now we state a proposition that finds a lower density for the set of $k$-layered numbers by having the smallest $k$-layered number.

\begin{proposition}\label{lowerd}
Let $k$ be a positive integer, and $n$ be the smallest $k$-layered number with prime factorization $p^{\alpha_1}_1p^{\alpha_2}_2\cdots p^{\alpha_m}_m$, where $p_1<p_2<\cdots<p_m$. Then $(L_{max}(p_1, p_2,\dots,p_m)+1)n$ is an upper bound for the difference of two consecutive $k$-layered numbers.
\end{proposition} 
\begin{proof}
Let $a$ and $b$ be two consecutive $k$-layered numbers. By the division algorithm, there exist non-negative integers $s\neq0$ and $r$ such that $a=sn+r$ and $0\leq r<s$. By the definition of $L_{max}$ and Proposition \ref{nw}, we know that there exists an integer $\ell$ with $s+1 \leq \ell \leq s+1+L_{max}(p_1,p_2,\dots,p_m)$ such that $\ell n$ is $k$-layered. Then $b-a\leq \ell n-sn\leq (L_{max}(p_1,p_2,\dots, p_m)+1)n$.
\end{proof}

\begin{corollary}\label{goodgcd}
Let $n$ be the smallest $k$-layered number. If $n$ with prime factorization $p^{\alpha_1}_1p^{\alpha_2}_2\cdots p^{\alpha_m}_m$ is a good $k$-layered number, then $b-a\leq(L_{max}(p_2,p_3,\dots,p_m)+1)n.$
\end{corollary}

By Propositions \ref{upper for L} and \ref{lowerd}, we have the following corollary.
  \begin{corollary}\label{lower density}
  Let $k$ be a positive integer, and $n$ be the smallest $k$-layered number. Let $a$ and $b$ be two consecutive $k$-layered numbers. Then, $$|b-a|\leq (p_1p_2\cdots p_m-(p_1-1)(p_2-1)\cdots(p_m-1)+1)n.$$
  \end{corollary}
  \begin{remark}
  Let $k1$ be a positive integer. Corollary \ref{lower density} gives a lower density for the set of $k$-layered numbers by having the smallest $k$-layered number.
  \end{remark}
Let $k$ be a positive integer such that $2\leq k \leq 5$, and $n$ be the smallest $k$-layered number. In Section \ref{the smallest}, we showed that $n$ possesses a prime factorization $p^{\alpha_1}_1p^{\alpha_2}_2\cdots p^{\alpha_m}_m$, where $p_i$ is the $i$th prime number for integers $1\leq i \leq m$. Then, by Corollaries \ref{good corollary} and \ref{goodgcd}, we conclude that $(L_{max}(p_2,\dots,p_m)+1)n$ is an upper bound for the difference of two consecutive $k$-layered number. Thus, by Proposition \ref{L at some value}, we have the following table.
\begin{table}[h!]
	\centering
    \label{tab:2}
    \begin{tabular}{|C{1.5cm}|C{4.5cm}|C{5cm}|}
    \hline
    \multicolumn{1}{|C{2cm}|}{ $k$ }
    & The smallest $k$-layered number & An upper bound for the difference of two consecutive $k$-layered numbers\\ \hline
    $1$ & $1$ & $0$ \\ \hline
    $2$ & $6$ & $12$ \\ \hline
    $3$ & $120$ & $360$ \\ \hline
    $4$ & $27720$ & $194040$ \\ \hline
    $5$ &$147026880$ & $1911349440$\\ \hline
   \end{tabular}
   \caption{An upper bound for the difference of two consecutive $k$-layered for $1\leq k\leq 5$.}
\label{tab:3}
   \end{table}
   \begin{remark}
    The numbers $282$ and $294$ are two consecutive Zumkeller numbers. Thus, by Table \ref{tab:2}, the number $12$ is a sharp upper bound for difference of two consecutive Zumkeller numbers.
   \end{remark}
 \section*{Acknowledgements}
I thank Pankaj Jyoti Mahanta and Prof. Manjil P. Saikia for reading the whole of this paper and pointing out some mirror mistakes, and also for very important suggestions for improving the structure of the paper. I also thank Prof. Michael Filaseta, Prof. Steven J. Miller, and Roohallah Mahkam for their comments that helped me improve the paper.

 \bigskip
\hrule
\bigskip

\noindent 2020 {\it Mathematics Subject Classification}. 11A25; 11D99; 11Y99.

\noindent \emph{Keywords: } $k$-layered numbers, $k$-multiperfect numbers, Near-perfect numbers, Perfect numbers, Practical numbers, Superabundant numbers, Zumkeller numbers.
\end{document}